\documentclass[notitlepage,11pt,draft]{article}
\usepackage{amssymb,amsfonts,amsmath,geometry,esint,comment,upgreek,cases, fancyhdr}

\catcode`\@=11
\@addtoreset{equation}{section}

\catcode`\@=12

\allowdisplaybreaks

\usepackage{color}

\def\R{\mathbb{R}}
\def\C{\mathbb{C}}

\def\Z{\mathbb{Z}}
\def\H{\mathbb H}
\def\w{\textsc U}

\newcommand{\overbar}[1]{\mkern 1mu\overline{\mkern-0.0mu#1\mkern-1mu}\mkern 1mu}
\newcommand{\ovR}{\overbar{\R}\!\,}

\def\f{\varphi}

\def\irn{\int\limits_{\R^2}}

\def\eps{\varepsilon}
\def\div{{\rm div}}
\def\e{\varepsilon}
\def\P{\mathcal{P}}
\def\F{\mathcal{F}}
\def\L{\mathcal{L}}
\def\S{\mathbb S}

\def\tange{T_{\!\w\!}Z}

\def\proof{\noindent{\textbf{Proof. }}}
\def\QED{\hfill {$\square$}\goodbreak \medskip}

\newtheorem{Theorem}{Theorem}[section]
\newtheorem{Lemma}[Theorem]{Lemma}

\newtheorem{Remark}[Theorem]{Remark}

\usepackage{mathtools}

\DeclareMathOperator{\artanh}{artanh}

\begin{document}
\title{Bubbles with constant mean curvature,  \\
and almost constant mean curvature, in the hyperbolic space}

\author{Gabriele Cora
 \and Roberta Musina}

\date{}

\maketitle

\noindent
{\footnotesize {\bf Abstract.}
 \it
Given a constant $k>1$, let $Z$ be the family of round spheres of radius $\artanh(k^{-1})$ in the hyperbolic space $\H^3$, so that 
any sphere in $Z$ has  mean curvature $k$.  We prove a crucial nondegeneracy result involving the manifold $Z$.
As an application, we provide sufficient conditions on a prescribed  function $\phi$ on $\H^3$, which
ensure the existence of a ${\cal C}^1$-curve, parametrized by $\eps\approx 0$,
of  embedded spheres in $\H^3$ having mean curvature $k +\eps\phi$
at each point. 
}



\normalsize

\bigskip

\section{Introduction}

Let $K$ be a given function on the hyperbolic space $\H^3$. 
The  $K$-bubble problem consists in finding a {\em $K$-bubble}, which is an immersed surface $u:\S^2\to \H^3$
having mean curvature $K$ at each point. 
Besides its independent interest, the significance of the $K$-bubble problem is  due to 
its  connection with the Plateau problem
for
disk-type parametric surfaces having prescribed mean curvature $K$ and contour $\Gamma$, see for instance
\cite{AR, Gu}. In the Euclidean case, the impact of $K$-bubbles on nonexistence and lack of compacteness phenomena
in the Plateau problem has been investigated in \cite{BrCo2, CMarma, CM_JFA}.

To look for $K$-bubbles in the hyperbolic setting  one can model $\H^3$ via the upper half-space $(\R^3_+, p_3^{-2}\delta_{h j})$
and  consider the elliptic system
\begin{equation}
\label{eq:bollehyp}
\Delta u - 2 u_3^{-1} {{G}(\nabla u)} =2 u_3^{-1}K(u) ~\!\partial_x u\wedge \partial_y u
\end{equation}
for functions $u=(u_1,u_2,u_3)\in {\cal C}^2(\S^2,\H^3)$. 
Here we used the stereographic projection to introduce local coordinates on $\S^2\equiv \R^2\cup\{\infty\}$ 
and put
\begin{equation}\label{eq:Gamma}
{G}_{\!\ell}(\nabla u):=\nabla u_3\cdot\nabla u_\ell- \frac{1}{2}|\nabla u|^2 \delta_{\ell 3}
=-\frac12u_3\sum_{h,j=1}^3{\Gamma}^\ell_{hj}(u)\nabla u_h\cdot\nabla u_j\,,\quad \ell=1,2,3\,,
\end{equation}
where ${\Gamma}^\ell_{hj}$ are the Christoffel symbols.
Any nonconstant solution $u$  to (\ref{eq:bollehyp}) is a {\em generalized} $K$-bubble in $\H^3$
(see Lemma \ref{L:conformal} in the Appendix and  \cite[Chapter 2]{GHR}),
that is, $u$ is
a conformal parametrization of a surface having  mean curvature $K(u)$, apart from
a finite number of branch points.  Once found a solution to (\ref{eq:bollehyp}), the next step should concern the study of the geometric regularity of the surface $u$, which might have self-intersections and branch points.

\smallskip

A remarkable feature of (\ref{eq:bollehyp}) is its variational structure, which means that
its  solutions are critical points of a certain energy functional $E$, see the Appendix for details. Because of their underlying geometrical meaning, both (\ref{eq:bollehyp}) and $E$ are invariant 
with respect to the action of 
M\"obius transformations. This produces
some lack of compactness phenomena, similar to those observed in the 
largely studied
 $K$-bubble problem, raised by S.T. Yau \cite{Yau82}, for surfaces in $\R^3$  (see for instance 
\cite{CM, CMna, Fe, TW} and references therein; see also \cite{An, CC, Yau} for related problems).
However, the hyperbolic $K$-bubble problem is definitively more challenging, due to the homogeneity properties that characterize the hyperbolic-area and the hyperbolic-volume functionals. 

The main differences between the Euclidean and the hyperbolic case are already evident 
when the prescribed curvature is a constant $k >0$ (the case $k<0$ is recovered by a change of orientation).
Any round sphere of radius $1/k$ in $\R^3$  can be parameterized by an embedded $k$-bubble,
which minimizes the energy functional 
$$
E_{\text{Eucl}}(u)=\frac12\irn  |\nabla u|^2~\! dz+\frac{2k}{3} \irn u\cdot \partial_x u\wedge\partial_yu~\!dz
$$
on the Nehari manifold $\{\, u\neq \text{const.}\ |\ E_{\text{Eucl}}'(u)u=0\,\}$, see \cite[Remark 2.6]{CM}. 
In contrast, no immersed hyperbolic $k$-bubble exists if $k\in(0,1]$, see for instance \cite[Theorem 10.1.3]{Lop}.
If $k>1$, then
any sphere in $\H^3$ of radius
$$
\rho_{k}: =\artanh \frac1{k}=\frac12\ln\frac{k+1}{k-1}
$$
can be parameterized by an embedding $U:\S^2\to\H^3$, which solves
\begin{equation*}\tag{$\mathcal P_0$}
\label{eq:problem}
\Delta u - 2 u_3^{-1} {{G}(\nabla u)} =2 u_3^{-1}k ~\!\partial_x u\wedge \partial_y u\quad \text{on $\R^2$}\,,
\end{equation*}
and which is a critical point of the energy functional
\begin{equation}
\label{eq:ene}
E_{\text{hyp}}(u)=\frac12\irn  u_3^{-2}|\nabla u|^2~\! dz-k \irn u_3^{-2}e_3\cdot \partial_x u\wedge\partial_yu~\!dz\, .
\end{equation}
As in the Euclidean case, the functional $E_{\text{hyp}}$ is unbounded from below (see Remark \ref{R:unbounded}). Thus $U$ does not minimize 
the energy $E_{\text{hyp}}$ on the Nehari manifold, which in fact fills $\{\,u\neq \text{const.}\,\}$.

\medskip
Besides their invariance with respect to M\"obius transformations, the system (\ref{eq:problem}) and the related 
energy $E_{\text{hyp}}$ are invariant  with respect to the 
$3$-dimensional group of hyperbolic translations as well.
Thus, any $k$-bubble generates a smooth $9$-dimensional manifold $Z$ 
of solutions
to (\ref{eq:problem}).  We 
explicitly describe the tangent space $T_ {\!U\!}Z$ at   $U\in Z$ in formula (\ref{eq:TUZ}).

As a further consequence of the invariances of problem (\ref{eq:problem}), any tangent direction
$\f\in T_{\!U\!}Z$ solves the elliptic system
\begin{equation}\label{eq:problem_lin}
\Delta {\f}-2U_3^{-1}{G}'(\nabla U)\nabla \f=-U_3^{-1}\f_3\Delta U
+2 {U}_3^{-1} k \big(\partial_x {\f}\wedge \partial_y {U}+\partial_x {U}\wedge \partial_y {\f}\big)\, , 
\end{equation}
which is obtained by linearizing (\ref{eq:problem}) at $U$. 

The next one is the main result of the present paper.

\begin{Theorem}[Nondegeneracy]
\label{T:main}
Let $U\in Z$. If $\f\in {\cal C}^2(\S^2,\R^3)$ solves the linear system (\ref{eq:problem_lin}), then $\f\in T_{\!U\!}Z$.
\end{Theorem}

In the Euclidean case the nondegeneracy of $k$-bubbles has been proved in \cite[Proposition 3.1]{Mu}. The 
proof of Theorem \ref{T:main} (see Section \ref{S:nondeg}),  is considerably more involved. It requires
the choice of a suitable orthogonal frame for functions in ${\cal C}^2(\S^2,\R^3)$ and the complete classification of solutions of two systems of linear elliptic differential equations, which might have an independent geometrical interest (see Lemmata \ref{L:tangsol}, 
\ref{L:omegasol}).

\medskip

As an application of Theorem \ref{T:main}, we provide sufficient conditions on 
a prescribed smooth function $\phi:\H^3\to \R$ that ensure the existence of embedded surfaces $\S^2\to\H^3$
having nonconstant mean curvature $k +\eps\phi$. 
Our existence results involve the notion of {\em stable critical point} already used in \cite{MuZu} and inspired from
\cite[Chapter 2]{AM} (see Subsection \ref{SS:stable}). The main tool
is a Lyapunov-Schmidt reduction technique  combined 
with variational arguments, in the spirit of  \cite{AM}. 

\begin{Theorem}
\label{T:existence}
Let $k>1$ and $\phi\in {\cal C}^1(\H^3)$ be given. Assume  that the function
\begin{equation}
\label{eq:Mel}
F^{\phi}_k(q):=\int\limits_{B^\H_{\rho_{k}}\!(q)}\phi(p)~\!d\H^3_p~,\quad F^{\phi}_k:\H^3\to\R
\end{equation}
has a stable 
critical point in an open set $A\Subset \H^3$. 
{For every $\eps\in\R$ close enough to $0$ there exist a point $q^\eps\in A$,
a conformal parametrization $\w_{\!q^\eps}$ of a sphere
of radius $\rho_k$ about $q^\eps$, and a conformally embedded $(k+\eps \phi)$-bubble $u^\eps$, such that
$\|u^\eps-\w_{\!q^\eps}\|_{{\cal C}^2}=O(\eps)$ as $\eps\to 0$.}

Moreover, any sequence $\eps_h\to 0$ has a subsequence $\eps_{h_j}$ such that 
$q^{\eps_{h_j}}$ converges to a critical point for $F^{\phi}_k$.
In particular, if 
$\hat q\in \Omega$ is the unique critical point for $F^{\phi}_k$ in $\overline \Omega$, then  
$u^\eps\to \w_{\!\hat q}$ 
in ${\cal C}^2(\ovR^2,\H^3)$. 

\end{Theorem}

\begin{Theorem}
\label{T:main2}
Assume that  $\phi\in {\cal C}^1(\H^3)$ has a stable 
critical point in an open set $A\Subset \H^3$. Then
there exists $k_0>1$ such that for any  $k>k_0$ 
and for every $\eps$ close enough to $0$, {there exists a conformally embedded  $(k+\eps \phi)$-bubble.}
\end{Theorem}

In Section  \ref{S:fdr} we first show that the existence of a critical point for $F^{\phi}_k(q)$ is a necessary condition
in Theorem \ref{T:existence}. Then we perform the dimension reduction and prove 
Theorems \ref{T:existence}, \ref{T:main2}. 
With respect to correspondent Euclidean results in \cite{CM}, a different choice of the functional setting allows us to weaker the regularity assumption on $\phi$ (from ${\cal C}^2$ to ${\cal C}^1$).

We conclude the paper with an Appendix in which we 
collect some partially known results about the variational approach to (\ref{eq:bollehyp}) and
prove a nonexistence result for (\ref{eq:bollehyp}) which, in particular, justifies the assumption on the existence of a critical point for $\phi$ in Theorem \ref{T:main2} .

\section{Notation and preliminaries}

The vector space $\R^n$ is endowed with the Euclidean scalar product $\xi\cdot \xi'$ and norm $|\xi|$.
We denote by $\{e_1, e_2, e_3\}$ the canonical basis
and by $\wedge$ the exterior product in $\R^3$.

We will often identify the complex number $z=x+iy$ with the vector $z=(x,y)\in \R^2$. Thus, 
$iz\equiv (-y,x)$, $z^2\equiv (x^2-y^2,2xy)$  and $z^{-1}\equiv|z|^{-2}(x,-y)$ if $z\neq 0$.

\medskip

In local coordinates induced by the 
stereographic projections from the north and the south poles, the round metric on the 
sphere $\S^2$ is given by $g_{hj}=\mu^2\delta_{hj}$, $d\S^2=\mu^2\!dz$, where
$$
\mu(z)=\frac{2}{1+|z|^2}\, .
$$
We  identify the compactified plane $\ovR^2$ with the sphere $\S^2$ through  
the inverse of the stereographic projection from the north pole, which is explicitly given by 
\begin{equation}
\label{eq:omega}
\omega(x,y)=(\mu x,\mu y, 1-\mu)\, .
\end{equation}
The identity $|\omega|^2\equiv 1$ trivially gives
$\omega\cdot\partial_x\omega\equiv 0$, $\omega\cdot\partial_y\omega\equiv 0$. 
We also notice that 
$\omega$ is a conformal (inward-pointing)
parametrization of the unit  sphere and satisfies
\begin{equation}
\label{eq:omega1}
\begin{cases}
\Delta \omega = 2 {\partial_{x}\omega} \wedge {\partial_{y}\omega}~,\quad -\Delta\omega=2\mu^2\omega\\
{\partial_{x}\omega} \cdot {\partial_{y}\omega}=0\,\\
|{\partial_{x}\omega}|^2 = |{\partial_{y}\omega}|^2 =\frac12|\nabla \omega|^2 = \mu^2~\!.
\\
\partial_{x}\omega \wedge \omega = \partial_{y} \omega~,\quad
\omega \wedge \partial_{y} \omega = \partial_{x}\omega~,\quad
\partial_{x}\omega \wedge \partial_y\omega = -\mu^2 \omega.
\end{cases}
\end{equation}

\subsection{The Poincar\'e half-space model}\label{SS:model}
We adopt
as model for the three dimensional hyperbolic space $\H^3$ the upper half-space
$\R^3_+=\{ (p_1,p_2,p_3)\in\R^3\ | \ p_3 > 0 \}$
endowed with the Riemannian metric $g_{h j}(z)= p_3^{-2}\delta_{h j}$\,. 

The hyperbolic distance $d_\H(p,q)$ in $\H^3$ is related to the Euclidean one  by 
$$\cosh d_\H(p,q) = 1 + \frac{|p-q|^2}{2 p_3 q_3}~\!,$$ 
and the hyperbolic ball $B^\H_\rho({p})$ centered at $p=(p_1, p_2,p_3)$ 
is the Euclidean ball of center 
$(p_1, p_2, p_3 \cosh \rho)$ and radius $p_3 \sinh \rho$.

\smallskip

If $F:\H^3\to\R$ is a differentiable function, then 
$\nabla^{\!\H} F(p)=p_3^2\nabla F(p)$, where $\nabla^{\!\H}$, $\nabla$
are the hyperbolic and the Euclidean gradients, respectively. In particular, 
$\nabla^{\!\H} F(p)=0$ if and only if $\nabla F(p)=0$.
The hyperbolic
volume form is related to the Euclidean one by $d{\H_p^3}=p_3^{-3}dp$.

\subsection{Stable critical points}
\label{SS:stable}

Let ${X}\in {\cal C}^1(\H^3)$ and let $\Omega\Subset \H^3$ be open. 
We say that ${X}$ has a stable critical point in $\Omega$ if
there exists $r>0$ such that any  function ${G}\in {\cal C}^1(\overline \Omega)$
satisfying 
$\displaystyle{\|{G}-{X}\|_{{\cal C}^1(\overline \Omega)}<r}$  has a critical point in $\Omega$.

As noticed in \cite{MuZu}, conditions to have the existence of a stable critical point $p\in \Omega$ for ${X}$
are easily given via elementary calculus. 
For instance, one can use Browder's topological degree theory or can assume that 
$$\displaystyle{{\min_{\partial \Omega} {X}>\min_{\Omega} {X}}}\quad\text{or}\quad
\displaystyle{\max_{\partial \Omega} {X}<\max_{\Omega} {X}}.
$$
Finally, if ${X}$ is of class ${\cal C}^2$ and $\Omega$ contains a 
nondegenerate critical point $p_0$ (i.e. the Hessian matrix of ${X}$ at $p_0$ is invertible),
then  $p_0$ is stable.

\subsection{Function spaces}
\label{SS:FS}

Any function $f$ on $\ovR^2$ is identified with $f\circ\omega^{-1}$, which is a function on $\S^2$.  
If no confusion can arise, from now on we write $f$ instead of $f\circ \omega^{-1}$. 

The Hilbertian norm on $L^2(\ovR^2,\R^n)\equiv L^2(\S^2,\R^n)$ is
given by 
$$
\|f\|_{2}^2=\int_{\R^2} |f|^2\,\mu^2\!dz<\infty~\!.
$$
Let $m\ge 0$. We endow 
$$
{\cal C}^m(\ovR^2,\R^n)=\{\,u\in {\cal C}^m(\R^2,\R^n)~|~ u(z^{-1})\in {\cal C}^m(\R^2,\R^n)\,\}\equiv {\cal C}^m(\S^2,\R^n)
$$ 
with the standard Banach space structure (we agree that
${\cal C}^m(\ovR^2,\R^n)={\cal C}^{\lfloor m\rfloor, m-\lfloor m\rfloor}(\ovR^2,\R^n)$ if $m$ is not an integer).
If $m$ is an integer, a norm in ${\cal C}^m(\ovR^2,\R^n)$ is given by
\begin{equation}
\label{eq:Cm_norm}
\|u\|_{{\cal C}^m}= \|u\|_\infty+\|\mu^{-m}\nabla^m u\|_\infty~\!.
\end{equation}
Since we adopt the upper space model for $\H^3$, we are allowed to write
$$
{\cal C}^m(\ovR^2,\H^3)={\cal C}^m(\ovR^2,\R^3_+)=\{u\in {\cal C}^m(\ovR^2,\R^3)~|~ u_3>0\}~\!,
$$
so that ${\cal C}^m(\ovR^2,\H^3)$ is an open subset of ${\cal C}^m(\ovR^2,\R^3)$.

\medskip

If $\psi,\f \in {\cal C}^1(\ovR^2,\R^3)$ and $\tau\in \R^2$, we put
$$
\nabla \psi \cdot \nabla \f = \partial_{x} \psi \cdot \partial_{x} \f + \partial_{y} \psi \cdot \partial_{y}\f~,\qquad
\tau\nabla\f= \tau_1\partial_{x}\f+\tau_2\partial_{y}\f
$$
(notice that $\tau\nabla \f(z)=d\f(z)\tau$ for any $z\in \R^2$). 
For instance, we have
$$
z^h\nabla\f=\begin{cases} \partial_{x}\f&\text{if $h=0$}\\ 
x {\partial_{x}\f} + y {\partial_{y}\f}&\text{if $h=1$}
\end{cases}~,\qquad
iz^h\nabla\f=\begin{cases} \partial_y\f&\text{if $h=0$}\\ 
-y {\partial_{x}\f} + x {\partial_{y}\f}&\text{if $h=1$}\, .
\end{cases}
$$
For future convenience we notice, without proof, that the next identities hold,
\begin{equation}
\label{eq:omega3}
\begin{cases}
~\!~\!~\!\partial_{x}\omega = e_1 - \omega_1 \omega - e_2 \wedge \omega\\
~\!~\!z\nabla\omega = e_3 - \omega_3 \omega\\
z^2\nabla \omega=  - (e_1 - \omega_1\omega + e_2 \wedge \omega)
\end{cases}
\quad\quad
\begin{cases}
~~~\partial_{y}\omega = e_2 - \omega_2 \omega + e_1 \wedge \omega\\
~\!~\!iz\nabla\omega= e_3 \wedge \omega,\\
iz^2\nabla \omega = e_2 -\omega_2 \omega - e_1 \wedge \omega~\!.
\end{cases}
\end{equation}

\medskip

The monograph \cite{Au} is our reference text for the theory of Sobolev spaces on Riemannian manifolds. 
In view of our purposes, it is important to notice that
$$H^1(\ovR^2,\R^n)= \{\, u\in H^1_{\rm loc}(\R^2,\R^n)~|~|\nabla u|+|u|~\mu\in L^2(\R^2)\,\}
\equiv H^1(\S^2,\R^n)\, .
$$

\medskip
We simply write $L^2(\ovR^2)$, ${\cal C}^m(\ovR^2)$ and $H^1(\ovR^2)$  instead of $L^2(\ovR^2,\R)$, ${\cal C}^m(\ovR^2,\R)$ and $H^1(\ovR^2,\R)$,
respectively.

\subsection{M\"obius transformations and hyperbolic translations}
\label{SS:mobius}
Transformations in $PGL(2,\C)$ are obtained by composing translations,
dilations,
rotations and complex inversion. Its Lie algebra admits as a basis the transforms 
$$
z\mapsto 1~,~~z\mapsto i~,~~ z\mapsto z~,~~z\mapsto iz~,~~z\mapsto z^2~,~~z\mapsto iz^2\,.
$$
Therefore, for any $u\in {\cal C}^1(\ovR^2,\H^3)$, the functions
$$
z^h\nabla u~,\quad iz^h\nabla u~,\qquad h=0,1,2\, , 
$$
span the tangent space to the  manifold $\{\, u\circ g~|~g\in PGL(2,\C)\, \}$ at $u$.

\medskip

Hyperbolic translations are obtained by composing a horizontal (Euclidean) translation
$p\mapsto p+ae_1+be_2$,  $a,b\in\R$
with an Euclidean homothety $p\mapsto tp$, $t>0$. 
Therefore, for any $u\in {\cal C}^m(\ovR^2,\H^3)$, the functions
$$
e_1~,\quad e_2~,\quad  u\,,
$$
span the tangent space to the  manifold $\{\, u_q~|~q\in\H^3\, \}$ at $u$, where
\begin{equation}
\label{eq:uq}
 u_q:=q_3 u+q-(q\cdot e_3)e_3\,.
\end{equation}

\section{Nondegeneracy of hyperbolic $k$-bubbles}
\label{S:nondeg}

The proof of Theorem \ref{T:main} needs some preliminary work.  
We put
\[
\w = r_{\!k}(\omega+ ke_3)\,,\quad r_{\!k}:=\sinh\rho_k=\frac1k\cosh\rho_{k}=\frac{1}{\sqrt{k^2-1}}\,,
\]
where $\omega$ is given by (\ref{eq:omega}). Since $\w$ is a conformal parametrization of the Euclidean sphere of radius $r_{\!k}$ about $kr_{\!k}e_3$, which coincides
with the hyperbolic sphere of radius $\rho_k$ about $e_3$, then $\w$ has  curvature 
$k$ and in fact it solves (\ref{eq:problem}). 
Accordingly with (\ref{eq:uq}), we put
\begin{equation}
\label{eq:wq}
\w_{\!q}:=q_3\w+q-(q\cdot e_3)~\!e_3
\end{equation}
(notice that $\w_{\!e_3}=\w$), and introduce the $9$-dimensional manifold
\begin{equation}\label{eq:Z}
Z=\big\{\, \w_{\!q}\circ g~|~g\in PGL(2,\C)~,~q\in\H^3\,\big\}~\!.
\end{equation}

\begin{Remark}
\label{R:unique}
Any surface $U\in Z$ is an embedding and solves {\em (\ref{eq:problem})}. 
Conversely, let
$U\in{\cal C}^2(\ovR^2,\H^3)$ be an embedding. If $U$ solves {\em (\ref{eq:problem})}, then
it is a $k$-bubble by Lemma \ref{L:conformal} and, thanks to an Alexandrov' type argument (see for instance \cite[Corollary 10.3.2]{Lop}) it parametrizes a  sphere of hyperbolic radius $\rho_k$ and Euclidean radius $r_k$. 
Since $U$ is conformal,  then 
$\Delta U=2r_{\!k}^{-1}\partial_xU\wedge \partial_yU$. Therefore
$U\in Z$ by the uniqueness result in \cite{BrCo2}.
\end{Remark}

By the remarks in Subsection \ref{SS:mobius}  and since $\nabla \w_{\!q}$ is proportional to 
$\nabla \omega$, we have that $T_{\!\w_{\!q}\!}Z=\tange$
for any $q\in\H^3$, and
\begin{equation}\label{eq:tangbasis}
\tange=\big\langle \{z^h\nabla \omega~, ~iz^h\nabla \omega~|~h=0,1,2\}\big\rangle\oplus\big\langle e_1, e_2, \w~\!\big\rangle\, .
\end{equation}
Moreover,  any tangent direction $\tau\in \tange$ solves (\ref{eq:problem_lin}).

It is convenient to split ${\cal C}^m(\ovR^2, \R^3)$ in the direct sum of its closed subspaces
\[
\langle\omega\rangle_{{\cal C}^m}^\perp 
: = \{\, \f \in {\cal C}^m(\ovR^2, \R^3) \ |\ \f \cdot \omega\equiv 0 \text{ on }\R^2\,\}\,,\quad
\langle\omega\rangle_{{\cal C}^m}:=\{\,\eta\omega~|~\eta\in {\cal C}^m(\ovR^2)\, \}~\!.
\]
Since
$\tange= \big(\tange\cap\langle\omega\rangle^\perp_{{C}^2}\big)\oplus \big(\tange\cap\langle\omega\rangle_{{C}^2}\big)$, from
(\ref{eq:omega3}) we infer another useful description of the tangent space, that is
\begin{equation}\label{eq:TUZ}
\tange=\big\{\, s-(s\cdot \omega)\omega+t\wedge\omega~|~s,t \in\R^3\,\big\}\oplus
\big\{\,(\alpha\cdot(k\omega+e_3))~\!\omega~|~\alpha\in\R^3\,\big\}\, .
\end{equation}

We now introduce
the differential operator
$$
J_0(u)= -\div(u_3^{-2}\nabla u)-u_3^{-3}|\nabla u|^2e_3+2ku_3^{-3}{\partial_{x}u}\wedge {\partial_{y}u}~\!.
$$
Notice that  $Z\subset \{J_0=0\}$. Further, let $I(z)=z^{-1}$. Since $J_0(u\circ I)=|z|^{-4}J_0(u)\circ I$ for any
$u\in {\cal C}^{2+m}(\ovR^2,\R^3)$, $m\ge 0$, then  $J_0$ is a ${\cal C}^1$ map
$$
J_0: {\cal C}^{2+m}(\ovR^2,\H^3)\to {\cal C}^m(\ovR^2,\R^3)~\!.
$$
We denote by $J'_0(u) : {\cal C}^{2+m}(\ovR^2,\R^3)\to {\cal C}^m(\ovR^2,\R^3)$
its differential at $u$.

Finally, $J_0(\w_{\!q}\circ g)=0$ for any $g\in PGL(2,\C)$, $q\in\H^3$, that implies $\tange \subseteq \ker J'_0(\w)$.
In order to prove Theorem \ref{T:main} it suffices to show that 
\[
\ker {J}'_0(\w)\subseteq \tange\, .
\]

\paragraph{Main computations}~\\
Recall that $\w=r_{\!k}(\omega+ke_3)$ solves 
$J_0(\w)=0$  to check that
\begin{multline*}
J'_0(\w)\varphi=
- \text{div}\big(\w_3^{-2}{\nabla \varphi}\big)\\+
2\w_3^{-3}\big[G'(\nabla\w)\nabla\f-\nabla U_3\nabla\f-\frac12\f_3\Delta \w
+k(\partial_x \f\wedge {\partial_{y}\w} + {\partial_{x}\w} \wedge \partial_y \f)\big]\, ,
\end{multline*}
where ${G}$ is given by (\ref{eq:Gamma}). Since $\nabla\omega_3=-\nabla\mu=\mu^2z$, thanks to (\ref{eq:omega1}) we have
\begin{align}
\label{eq:Jsimpl}
&\begin{aligned}
r_{\!k}^2J'_0&(\w)\varphi=
- \text{div}\big((\omega_3+k)^{-2}{\nabla \varphi}\big)\\
 &+ 2(\omega_3 + k)^{-3}
\big[\big({G}'(\nabla \omega)\nabla\f -\mu^2z\nabla \f\big)\!+\mu^2 \varphi_3\omega +  {k} 
(\partial_x \f\wedge {\partial_{y}\omega} + {\partial_{x}\omega} \wedge \partial_y \f)\big]\, ,
\end{aligned}\\
\label{eq:Gsimpl}
&{G}'(\nabla \omega)\nabla\f -\mu^2z\nabla \f=\nabla\f_3\nabla\omega-(\nabla\f\cdot\nabla\omega)e_3\, .
\end{align}
To rewrite (\ref{eq:Jsimpl})  in a less obscure form, we decompose any $\f\in {\cal C}^m(\ovR^2,\R^3)$, $m\ge 0$, as
\begin{equation}\label{eq:omegadec}
\f ={\P\!\f} + (\f \cdot \omega) \omega\, , \qquad 
{\P\!\f} := \f - (\f \cdot \omega)\omega=\mu^{-2}\big((\f\cdot\partial_x\omega)\partial_x\omega+(\f\cdot \partial_y\omega)\partial_y\omega\big)~\!.
\end{equation}
Accordingly, for $\f\in {\cal C}^{2}(\ovR^2,\R^3)$ we have 
$$J'_0(\w)\f= \mathcal P\big(J_0'(\w)\f\big)+(J_0'(\w)\f\cdot\omega)\omega\, ,
$$
so that we can reconstruct $J'_0(\w)\f\in {\cal C}^0(\ovR^2,\R^3)$ by providing explicit expressions for
$\mathcal P\big(J_0'(\w)\f\big)$ and $J_0'(\w)\f\cdot\omega$, separately.  This will be done in the next Lemma.

\begin{Lemma}
\label{L:J_split}
Let $\f \in {\cal C}^2(\ovR^2, \R^3)$. Then
\begin{align}
\label{eq:J'dec1}
&r_{\!k}^2\mathcal P\big(J_0'(\w)\f\big) =
{\P} \Big(\!-\div\Big(\frac{\nabla {\P\!\f}}{(\omega_3 + k)^2} \Big)\Big)+ 
\frac{2\mu^2}{(\omega_3 + k)^3}~\!(iz\nabla  {\P\!\f}) \wedge \omega - \frac{2\mu^2}{(\omega_3 + k)^2} ~\!{\P\!\f}\, ,\\
\label{eq:J'dec2}
&r_{\!k}^2(J_0'(\w)\f)\cdot\omega= -\div\Big(\frac{\nabla (\f\!\cdot\! \omega)}{(\omega_3 + k)^2}\Big) - 
\frac{2k\, \mu^2}{(\omega_3 + k)^3}~\!(\f\!\cdot\! \omega)~\!.
\end{align}
\end{Lemma}

\proof
We introduce the differential operator 
$L=-\div\big((\omega_3 + k)^{-2}\nabla~\big)
$
and start to prove (\ref{eq:J'dec2}) by noticing that 
\begin{equation}\label{eq:Lproj}
L\f\cdot \omega=L(\f\cdot\omega)+2(\omega_3+k)^{-3}\big[(\omega_3+k)\nabla\f\cdot\nabla\omega-\mu^2\f\cdot(z\nabla\omega) - \mu^2(\omega_3 + k)(\f \cdot \omega)\big]\, .
\end{equation}
Recalling that $\omega$ is pointwise orthogonal to $\partial_x\omega, \partial_y\omega$,
from (\ref{eq:Gsimpl}) we obtain
$$
\big({G}'(\nabla \omega)\nabla\f-\mu^2z\nabla \f\big)\cdot\omega =
-(\nabla\f\cdot\nabla\omega)\omega_3~\!.
$$
Further, by (\ref{eq:omega1}) we have
$(\partial_x \f\wedge {\partial_{y}\omega} + {\partial_{x}\omega} \wedge \partial_y \f)\cdot\omega=
-\nabla\f\cdot\nabla\omega$. Finally, we obtain
$$
r_{\!k}^2(J_0'(\w)\f)\cdot\omega=L(\f\cdot\omega)-2(\omega_3+k)^{-3}\mu^2
\big[\f\cdot(z\nabla\omega)-\f_3 + (\omega_3 + k) (\f\cdot \omega)\big]\, ,
$$
and (\ref{eq:J'dec2}) follows, because $e_3=z\nabla\omega+\omega_3\omega$, see (\ref{eq:omega3}).

\medskip

Next, using the equivalent formulation
\[
\begin{aligned}
\w_3^2 J'_0(\w)\varphi =-\Delta \f+ 2(\omega_3 + k)^{-1}\left[G'(\nabla\omega)\nabla\varphi + \mu^2 \omega \varphi_3+ k(\partial_{x}\varphi \wedge {\partial_{y}\omega} + {\partial_{x}\omega} \wedge \partial_{y}\varphi)\right]
\end{aligned}
\]
we find that, for $\f = \eta \omega$, $\eta \in {\cal C}^2(\ovR^2)$, it holds
\[
\w_3^2 J'_0(\w)(\eta \omega)\cdot \partial_x \omega = \w_3^2 J'_0(\w)(\eta \omega)\cdot \partial_y \omega = 0\, ,
\] 
whence we infer
\begin{equation}\label{eq:omegaproj}
\P\big(J_0'(\w)(\f - \P\f)\big) = 0\, , \quad \text{for every }\f \in {\cal C}^2(\ovR^2, \R^3)\,.
\end{equation}

Thanks to (\ref{eq:J'dec2}) and \eqref{eq:omegaproj} we get $\P\big(J_0'(\w)\f\big)=J_0'(\w)(\P\f)$, thus to conclude the proof we can assume that $\f=\P\f$. Since
$\f$ is pointwise orthogonal to $\omega$, we trivially have
$$
\partial_x\f\cdot\omega=-\f\cdot\partial_x\omega\,,\quad \partial_y\f\cdot\omega=-\f\cdot\partial_y\omega\, .
$$
We start to handle (\ref{eq:Gsimpl}). From $e_3=z\nabla\omega+\omega_3\omega$ we get
\begin{multline*}
(G'(\nabla\omega)\nabla\f-\mu^2z\nabla\f)+\omega_3(\nabla\f\cdot\nabla\omega)\omega=
\nabla\f_3\nabla\omega- (\nabla\f\cdot\nabla\omega)z\nabla\omega\\
=
\big(\partial_x\f_3-x (\nabla\f\cdot\nabla\omega)\big)\partial_x\omega+\big(\partial_y\f_3-y (\nabla\f\cdot\nabla\omega)\big)\partial_y\omega\, .
\end{multline*}
Further,
\begin{multline*}
\partial_x\f_3-x (\nabla\f\cdot\nabla\omega)=\partial_x\f\cdot(z\nabla\omega+\omega_3\omega)-x(\nabla\f\cdot\nabla\omega)\\
=
\big(\partial_x\f\cdot(z\nabla\omega)-x (\nabla\f\cdot\nabla\omega)\big)-\omega_3\f\cdot\partial_x\omega
=-(iz\nabla\f)\cdot\partial_y\omega -\omega_3\f\cdot\partial_x\omega\, .
\end{multline*}
In a similar way one can check that $\partial_y\f_3-y(\nabla\f\cdot\nabla\omega)= (iz\nabla\f)\cdot\partial_x\omega -\omega_3\f\cdot\partial_y\omega$,
thus
\[
{G}'(\nabla \omega)\nabla\f -\mu^2  z\nabla\f =  
\mu^2 (iz \nabla\f) \wedge \omega -\omega_3(\nabla \f \cdot\nabla \omega)\omega - \mu^2\omega_3\f.
\]
Next, using  (\ref{eq:omega1}) we can compute
\begin{gather*}
\partial_x\f\wedge\partial_y\omega=\partial_x\f\wedge(\partial_x\omega\wedge\omega)=
-(\f\cdot\partial_x\omega)\partial_x\omega-(\partial_x\f\cdot\partial_x\omega)\omega\\
\partial_x\omega\wedge\partial_y\f=(\omega\wedge\partial_y\omega)\wedge \partial_y\f=
-(\f\cdot\partial_y\omega)\partial_y\omega-(\partial_y\f\cdot\partial_y\omega)\omega,
\end{gather*}
which give the identity
\begin{equation}
\label{eq:misto}
\partial_x\f\wedge\partial_y\omega+\partial_x\omega\wedge\partial_y\f=-\mu^2\f-(\nabla\f\cdot\nabla\omega)\omega~\!,
\end{equation}
that holds for any $\f\in \langle\omega\rangle_{{\cal C}^m}^\perp$. 

Putting together the above informations we arrive at
$$
r_{\!k}^2J_0'(\w)\f=L\f+
\frac{2\mu^2}{(\omega_3 + k)^3}~\!(iz\nabla  \f) \wedge \omega - \frac{2\mu^2}{(\omega_3 + k)^2} ~\!\f
+\frac{2}{(\omega_3+k)^3}\big[\mu^2\f_3-(\omega_3+k)\nabla\f\cdot\nabla\omega\big]\omega\, .
$$
Using \eqref{eq:Lproj} and $\f_3 = \f \cdot (z\nabla \omega)$, we conclude the proof. \QED

Thanks to Lemma \ref{L:J_split}  we can study the system $J'_0(\w)\f=0$ separately, on $\langle\omega\rangle_{{\cal C}^m}^\perp$ 
first, and on $\langle\omega\rangle_{{\cal C}^m}$ later.  In fact, $\f\in \ker J'_0(\w)$ if and only if the pair of functions
$$
\psi:= \P\!\f 
\in \langle\omega\rangle_{{\cal C}^2}^\perp\subset {\cal C}^2(\ovR^2,\R^3)~,\qquad \eta:=\f\!\cdot\!\omega\in  
{\cal C}^2(\ovR^2)\, , $$ 
solves 
\begin{subnumcases}
{\label{eq:both}}
\displaystyle{
\P\Big(-\div\Big(\frac{\nabla\psi}{(\omega_3 + k)^2} \Big)\Big)+ 
\frac{2\mu^2}{(\omega_3 + k)^3}~\!(iz\nabla \psi) \wedge \omega = \frac{2\mu^2}{(\omega_3 + k)^2} ~\!{\psi}}\, ,\label{eq:prima} 
\\ 
\displaystyle{-\div\Big(\frac{\nabla \eta}{(\omega_3 + k)^2}\Big) =
\frac{2k\, \mu^2 }{(\omega_3 + k)^3}~\!\eta}~\!.
\label{eq:seconda}
\end{subnumcases}

\medskip

We begin by facing problem (\ref{eq:prima}). Firstly, we show that the quadratic form associated to the differential operator $J'_0(\w)$ is nonnegative on $\langle\omega\rangle_{{\cal C}^2}^\perp$.

\begin{Lemma}
\label{L:semidefpos}
Let $\psi \in \langle\omega\rangle_{{\cal C}^2}^\perp$. Then
\[
\irn J_0'(\w)\psi\cdot\psi~\!dz = r_k^{-2} \int\limits_{\R^2} 
\frac{(\partial_ x \psi \cdot \partial_ x\omega - \partial_y\psi \cdot \partial_y\omega)^2 + (\partial_ x\psi\cdot \partial_y\omega + \partial_y\psi \cdot \partial_ x\omega)^2}{\mu^2(\omega_3 + k)^2} ~\!dz
~\!.
\]
\end{Lemma}
\proof
Since $J_0'(\w)\psi\cdot\psi=\P\big(J_0'(\w)\psi\big)\cdot\psi$ and $\P\psi=\psi$, formula (\ref{eq:J'dec1}) gives
$$
r_k^{2} \irn J_0'(\w)\psi\cdot\psi~\!dz =
\irn\frac{|\nabla\psi|^2}{(\omega_3+k)^2}~\!dz+
2\int\limits_{\R^2}\frac{\psi\cdot(iz \nabla  \psi) \wedge \omega}{(\omega_3 + k)^3}~\mu^2\!d z 
-2\int\limits_{\R^2}\frac{|\psi|^2 }{(\omega_3 + k)^2}~\mu^2\!d z.
$$
Now we prove the identity 
\begin{equation} \label{eq:integralid}
B_\psi:=2\int\limits_{\R^2}\frac{\psi\cdot(iz \nabla  \psi) \wedge \omega}{(\omega_3 + k)^3}~\mu^2\!d z = 2
\int\limits_{\R^2}\frac{\omega \cdot \partial_x \psi\wedge \partial_y \psi}{(\omega_3 + k)^2}\, d z +  
\int\limits_{\R^2}\frac{|\psi|^2 }{(\omega_3 + k)^2}~\mu^2\!d z\, .
\end{equation}
We use polar coordinates $\rho,\theta$ on $\R^2$ and notice that
$\partial_\theta\psi = iz \nabla \psi$. 
From $\rho\mu^2=\partial_\rho\omega_3$ we get
\[
\begin{aligned}
B_\psi=&-
\int\limits_0^{2\pi}d\theta
\int\limits_0^{\infty} ~\!(\psi \cdot \partial_\theta\psi \wedge \omega)\partial_\rho (\omega_3 + k)^{-2}\, d \rho \\
=&
\int\limits_0^{\infty}\!d\rho\! \int\limits_0^{2\pi}\frac{\omega \cdot \partial_\rho \psi \wedge \partial_\theta\psi  
 - \psi \cdot  \partial_\rho \omega\wedge \partial_\theta\psi }{(\omega_3 + k)^2}d\theta +
 \int\limits_0^{\infty}\!d\rho\! \int\limits_0^{2\pi}\frac{\partial_{\rho\theta}\psi  \cdot \omega \wedge \psi }{(\omega_3 + k)^2}d\theta\, 
\\
=&\int\limits_0^{\infty}\!d\rho\! \int\limits_0^{2\pi}\frac{\omega \cdot \partial_\rho \psi \wedge \partial_\theta\psi  
 - \psi \cdot  \partial_\rho \omega\wedge \partial_\theta\psi }{(\omega_3 + k)^2}d\theta+
 \int\limits_0^{\infty}\!d\rho\! \int\limits_0^{2\pi}\frac{  \omega \cdot \partial_\rho \psi\wedge \partial_{\theta} \psi  - \psi \cdot\partial_\rho \psi \wedge \partial_\theta\omega }{(\omega_3 + k)^2}d\theta~\!.
\end{aligned}
\]
Using the elementary identity $\partial_\rho\alpha \wedge \partial_\theta\beta=\rho(\partial_x\alpha \wedge \partial_y\beta)$, we see that 
\[
B_\psi = 
2\int\limits_{\R^2}\frac{\omega \cdot 
\partial_x \psi\wedge \partial_y \psi}{(\omega_3 + k)^2}\, d z - 
\int\limits_{\R^2}\frac{\psi \cdot(\partial_x \omega\wedge \partial_y\psi + \partial_x \psi\wedge \partial_y \omega)}{(\omega_3 + k)^2}\, d z\,,
\]
and  (\ref{eq:integralid}) follows from (\ref{eq:misto}) (with $\f$ replaced by $\psi$).

Thanks to (\ref{eq:integralid}), we have the identity
$$
r_k^{2} \irn J_0'(\w)\psi\cdot\psi~\!dz=\irn\frac{|\nabla\psi|^2+2\omega\cdot\partial_x\psi\wedge\partial_y\psi-\mu^2|\psi|^2}{
(\omega_3+k)^2}~\!dz,
$$
so that we only need to handle the function
$$
b_\psi:=|\nabla\psi|^2+2\omega\cdot\partial_x\psi\wedge\partial_y\psi-\mu^2|\psi|^2.
$$
We decompose
$\partial_x\psi$ and $\partial_y\psi$ accordingly with \eqref{eq:omegadec}, to obtain
$$
\begin{gathered}
\mu^2\partial_x\psi=(\partial_x \psi\cdot \partial_x\omega)\omega_x + (\partial_x\psi\cdot \partial_y\omega)\omega_y-
\mu^2(\psi\cdot\partial_x\omega)\omega\,,\\
\mu^2\partial_y\psi=(\partial_y \psi\cdot \partial_x\omega)\omega_x + (\partial_y\psi\cdot \partial_y\omega)\omega_y-
\mu^2(\psi\cdot\partial_y\omega)\omega\, ,
\end{gathered}
$$
respectively. Since
$|\nabla\psi|^2=|\partial_x\psi|^2+|\partial_y\psi|^2$, we infer 
\[
\mu^2\big(|\nabla \psi|^2-\mu^2|\psi|^2) =(\partial_x \psi\cdot \partial_x\omega)^2 + (\partial_x\psi\cdot \partial_y\omega)^2+ (\partial_y \psi\cdot \partial_x\omega)^2+ (\partial_y \psi\cdot\partial_y\omega)^2\,.
\]
Writing $\mu^2\omega=-\partial_x\omega\wedge\partial_y\omega$, see (\ref{eq:omega1}), we get
$$
\mu^2\omega\cdot (\partial_x \psi\wedge \partial_y \psi )
=-(\partial_x\psi\cdot \partial_x \omega)(\partial_y \psi\cdot \partial_y\omega) +(\partial_x \psi\cdot \partial_y\omega)( \partial_y\psi\cdot \partial_x\omega)~\!,
$$
from which it readily follows that
$\mu^2 b_\psi=(\partial_x\psi\cdot\partial_x\omega -  \partial_y\psi\cdot\partial_y\omega)^2 + (\partial_x\psi\cdot \partial_y \omega+ \partial_y\psi\cdot \partial_x\omega)^2$. The proof is complete.
\QED

\begin{Lemma}\label{L:tangsol}
Let $\psi \in {\cal C}^2(\ovR^2,\R^3)$ be a solution to (\ref{eq:prima}). There exist $s,t \in \R^3$
such that 
$$\psi= s - (s\cdot \omega) \omega + t \wedge \omega\, ,$$
and thus $\psi \in \tange\cap\langle\omega\rangle^\perp_{{C}^2}=\{\, s - (s\cdot \omega) \omega + t \wedge \omega \ |\ s, t \in \R^3\, \}$.
\end{Lemma}
\proof
From (\ref{eq:prima}) it immediately follows that $\psi$ is pointwise orthogonal to $\omega$, which implies $\psi\in \langle\omega\rangle_{{\cal C}^2}^\perp$. Since $\mathcal P\psi=\psi$, then 
$J'_0(\w)\psi=0$ by (\ref{eq:J'dec1}) and (\ref{eq:J'dec2}), hence
\begin{equation}\label{eq:passsyst}
\begin{cases}
\partial_x \psi \cdot \partial_x\omega - \partial_y\psi \cdot \partial_y\omega = 0\\
\partial_x \psi \cdot \partial_y \omega + \partial_y\psi \cdot \partial_x \omega = 0
\end{cases}
\end{equation}
by Lemma \ref{L:semidefpos}. Since $\psi\in\langle\partial_x\omega,\partial_y\omega\rangle$ pointwise on $\R^2$, we can write
$$
\psi=f\nabla\omega\, ,\quad \text{where}\quad 
f:=\mu^{-2}(\psi\cdot \partial_x\omega, \psi\cdot \partial_y\omega)\in {\cal C}^2(\ovR^2, \R^2)\, .
$$
We identify $f$ with a complex valued function.
A direct computation based on \eqref{eq:omega1} shows that $\psi$ solves \eqref{eq:passsyst} if and only if $f$ 
solves $\partial_x f +  i \partial_y f = 0$ on $\R^2$. In polar coordinates we have that 
\begin{equation}\label{eq:tangode}
\rho \partial_\rho f + i \partial_\theta f = 0.
\end{equation}
For every $\rho>0$ we expand the periodic function $f(\rho,\cdot)$ in Fourier series,
$$
\displaystyle f(\rho, \theta) = \sum_{h\in \Z} \gamma_h(\rho) e^{i h \theta }~,
\quad 
\gamma_h (\rho) = \frac{1}{2\pi}\int\limits_0^{2\pi} f(\rho, \theta)e^{- ih\theta}d \theta\, .
$$
The coefficients $\gamma_h$ are complex-valued functions on the half-line $\R_+$ that solve
\[
\gamma'_h - h \gamma_h = 0\, ,
\]
because of  \eqref{eq:tangode}. Thus for every $h \in \Z$ there exists $a_h\in\C$
such that 
$\gamma_h(\rho) = a_h\rho^h$.
Now recall that $\mu\psi\in L^2(\R^2,\R^3)$. Since
\[
\int\limits_{\R^2}\mu^2|\psi|^2\, dz = \int\limits_{\R^2}\mu^4|f|^2\, d z 
\ge 2\pi \int_0^\infty\mu^4\rho |\gamma_h|^2\, d\rho=
a_h^2\irn \mu^4  |z|^{2h}\, dz \, , 
\quad \forall\ \! h \in \Z\, ,
\]
we infer that $\gamma_h = 0$ for every $h \neq 0,1,2$. Thus
$f(z) = \sum\limits_{h=0}^2 a_h z^h$, that is
$\psi= \sum\limits_{h=0}^2 a_h z^h\nabla\omega$, 
and in particular the space of solutions of (\ref{eq:prima}) has (real) dimension $6$. The conclusion of the
proof follows from the relations (\ref{eq:omega3}). \QED

\begin{Lemma}\label{L:omegasol}
Let $\eta \in{\cal C}^2(\ovR^2)$ be a solution to (\ref{eq:seconda}). There exists $\alpha\in \R^3$
such that
$$\eta=  \alpha\cdot(k\omega+e_3)~\!
$$
and thus $\eta \omega \in \tange\cap\langle\omega\rangle_{\mathcal{C}^2}=\{\, (\alpha\cdot(k\omega+e_3))~\!\omega~|~\alpha\in\R^3\, \}$.
\end{Lemma}

\proof
First of all, we notice that $\alpha\cdot(k\omega+e_3)$ solves (\ref{eq:seconda}) for any $\alpha\in\R^3$.

By the Hilbert--Schmidt theorem, the eigenvalue problem
\begin{equation}
\label{eq:E}
- \div\Big(\frac{\nabla \eta}{(\omega_3 + k)^2}\Big) =\frac{ \lambda\, \mu^2 }{(\omega_3 + k)^3}~\!\eta\quad\text{on $\R^2$}~\!,\quad
\eta\in {\cal C}^2(\R^2)\, ,
\end{equation}
has a non decreasing, divergent sequence $(\lambda_h)_{h\ge 0}$ of eigenvalues which correspond to
critical levels of the quotient
$$
R(\eta):=\frac{\displaystyle{
\int_{\R^2} \frac{|\nabla \eta|^2}{(\omega_3 + k)^2}~\!dz
}}
{\displaystyle{\int_{\R^2}\frac{{|\eta|^2}}{(\omega_3 + k)^3}~\!\mu^2\!dz}}~,\quad \eta\in H^1(\ovR^2)\setminus\{0\}\, .
$$
Clearly, $\lambda_0=0$ is simple, and its 
eigenfunctions are constant functions. We claim that the next eigenvalue is $2k$, and that its eigenspace has dimension $3$, that concludes the proof.

To this goal, we use the functional change 
$$
{\eta(z)= \frac{\mu(z)}{\mu(c_k z)}\Phi(c_k z)}\,,\quad c_k:=e^{\rho_k}=\sqrt{\frac{k+1}{k-1}}~\!.
$$
By a direct computation involving the identity
$(\omega_3(z) + k )\mu(c_k z) =  (k-1)\mu(z)$ and integration by parts, one gets
$$
\lambda_1=\inf_{\eta\in {\cal C}^2(\ovR^2)\setminus\{0\}\atop \scriptstyle\int_{\R^2}\frac{\eta~ \mu^2\!dz}{(\omega_3 + k)^3}=0}
R(\eta)=
2k+ 
\inf_{\Phi\in {\cal C}^2(\ovR^2)\setminus\{0\} \atop \scriptstyle\int_{\R^2}{\Phi}\mu^2\!dz=0}
\frac{\displaystyle\int_{\R^2} |\nabla \Phi|^2 dz- 2\displaystyle \int_{\R^2}{|\Phi|^2}\mu^2\!dz}{\displaystyle\int_{\R^2}\frac{|\Phi|^2}{(k-\omega_3)}\, \mu^2\! dz}\,.
$$
On the other hand, it is well known that 
$$
\min_{\Phi\in {\cal C}^2(\ovR^2)\setminus\{0\} \atop  \int_{\R^2}{\Phi}\mu^2\!dz=0}
\frac{\displaystyle\int_{\R^2} |\nabla \Phi|^2 dz}{\displaystyle \int_{\R^2}{|\Phi|^2}\mu^2\!dz}=2
$$
is the first nontrivial eigenvalue for the Laplace-Beltrami operator on the sphere and that its eigenspace has
dimension $3$, see for instance \cite{Au}. This concludes the proof.
\QED

\begin{Remark}
\label{R:spectrum}
The third eigenvalue $\lambda_2$ of (\ref{eq:E})  verifies $\lambda_2 >2k$ by Lemma \ref{L:omegasol}, and 
$$
\lambda_2=\min\Big\{\,  R(\eta)~\Big|~~
\irn  \frac{{\eta}}{(\omega_3 + k)^3}~\!\mu^2\!dz=\irn  \frac{{\eta}(k\omega_j+\delta_{j3})}{(\omega_3 + k)^3}~\!\mu^2\!dz=0~,~~j=1,2,3
\, \Big\}\, .
$$
\end{Remark}

\noindent
{\bf Proof of Theorem \ref{T:main}}
In fact, we only have to sum up the argument. Let $U \in Z$. Thanks to \eqref{eq:Z}, $U = \w_{\!q} \circ g$ for some $q \in \H^3$, $g \in PGL(2,\C)$. Since
\[
T_{\!{\w_{\!q} \circ g}}Z = \tange\circ g\,, \quad \ker J'_0(\w_{\!q} \circ g) = \ker J_0'(\w)\circ g\, , \quad \text{for every }q \in \H^3\,, \ g \in PGL(2,\C)\, ,
\]
it suffices to consider the case $U = \w$. 

If $\f\in {\cal C}^2(\ovR^2,\R^3)$ solves (\ref{eq:problem_lin})
then $J'_0(\w)\f=0$, which means
$\mathcal P\left(J_0'(\w)\f\right)=0$ and $\left(J_0'(\w)\f\right)\cdot\omega=0$.
From Lemma \ref{L:J_split} we infer that  $\mathcal P\f$ solves (\ref{eq:prima}) and that $\f\cdot\omega$ solves
(\ref{eq:seconda}). Therefore, Lemmata \ref{L:tangsol}, \ref{L:omegasol} give the existence of $s,t,\alpha\in\R^3$
such that
$$
 \mathcal P\f=s - (s\cdot \omega) \omega + t \wedge \omega~,\quad \f\cdot\omega=\alpha\cdot(k\omega+e_3)\, .
$$
Thus $\f=\mathcal P\f+(\f\cdot\omega)\omega\in \tange$ by (\ref{eq:TUZ}), which concludes the proof.
\QED

\subsection{Further results on the operator $J'_0(\w)$}
\label{S:Fred}

To shorten notation we put
$$
H^1 = H^1(\ovR^2,\R^3)~\!.
$$
Since integration by parts gives
$$
\irn - \text{div}\Big(\frac{\nabla \varphi}{(\omega_3 + k)^{2}}\Big)\cdot\psi~\!dz=
\irn \frac{\nabla \varphi\cdot\nabla\psi}{(\omega_3 + k)^{2}}~\!dz\ , \quad \f,\psi\in {\cal C}^2(\ovR^2,\R^3)~\!,
$$
the quadratic form 
\begin{equation}
\label{eq:aggiunta}
(\f,\psi)\mapsto \irn J'_0(\w)\f\cdot\psi~\!dz
\end{equation}  
can be extended to a continuous bilinear form
$H^1\times H^1\to\R$
via a density argument. It can be checked by direct computations (see also Remark \ref{R:aggiunta}), that the quadratic form in (\ref{eq:aggiunta}) is self-adjoint on $H^1$, that is,
\begin{equation}
\label{eq:SA}
\irn J'_0(\w)\f\cdot\psi~\!dz=\irn J'_0(\w)\psi\cdot\f~\!dz\quad\text{for any $\f,\psi\in H^1$.}
\end{equation}

Since $\tange$ is a subspace of $L^2(\ovR^2, \R^3)\equiv L^2(\S^2,\R^3)$, we are allowed to put
$$
\tange^\perp := \Big\{\, f \in L^2(\ovR^2, \R^3) \ \Big|\ \int\limits_{\R^2} f\cdot \tau~\!\mu^2\!dz= 0\, , \ \forall\ \tau \in \tange\, \Big\}\, .
$$
To shorten notation we introduce on $L^2(\ovR^2, \R^3)$ the equivalent scalar product
$$
(f, \psi)_* = \int\limits_{\R^2} \frac{\P f \cdot\P\psi}{(\omega_3 + k)^2}\, \mu^2\!dz + \int\limits_{\R^2} 
\frac{(f \cdot\omega)(\psi\cdot \omega)}{(\omega_3 + k)^3}\, \mu^2\!dz
$$
 and the subspaces
\begin{align*}
\tange^\perp_* :=& \big\{\, f \in L^2(\ovR^2, \R^3) \ |\ (f, \tau)_* = 0, \ \forall\, \tau \in \tange\,  \big\}\, ,\\
N_* :=&\langle\omega\rangle^\perp_*=
 \big\{\, f \in L^2(\ovR^2, \R^3) \ |\ (f, \omega)_*   = 0\,  \big\}~\!.
\end{align*}

We are in position to state the main result of this section.

\begin{Lemma}\label{L:Fredholm}
Let $q\in \H^3$.
For any $v\in \tange^\perp$,
there exists $\f_v\in H^1\cap \tange^\perp_*\cap N_*$ such that 
\begin{equation}
\label{eq:inverter0}
J'_0(\w_{\!q})\f_v=v~\!\mu^2\qquad\text{on $\R^2$.}
\end{equation}
If in addition $v\in {\cal C}^m(\ovR^2,\R^3)$ for some $m\in(0,1)$, then 
$\f_v\in {\cal C}^{2+m}(\ovR^2,\R^3)$.
\end{Lemma}

In view of Lemma \ref{L:J_split}, we split the proof of Lemma \ref{L:Fredholm} in few steps. 

\begin{Lemma}\label{L:tangnonhom}
Let $v \in \tange^\perp$ be such that $v \cdot \omega \equiv 0$ on $\R^2$. There exists  
$\f \in H^1\cap \tange^\perp_*$ such that $\f\cdot \omega\equiv 0$ on $\R^2$ and
\begin{equation}
\label{eq:inverter0bis}
J'_0(\w)\f=v~\!\mu^2\qquad\text{on $\R^2$.}
\end{equation}
\end{Lemma}

\proof 
We introduce
\[
X:= \big\{\, \psi\in H^1~\big|~ \psi \cdot\omega\equiv 0~~\text{on $\R^2$}\, \big\}\cap \tange^\perp_*~\!,
\]
which is a closed subspace of $H^1$.
Notice that $\psi= \P\psi$ for any $\psi\in X$ and moreover
$$
\irn J'_0(\w)\psi\cdot\psi~\!dz= \int\limits_{\R^2}
\frac{|\nabla {\psi}|^2}{(\omega_3 + k)^2}~\!dz + 2
\int\limits_{\R^2} \Big(
\frac{(\psi\cdot iz\nabla  {\psi})\wedge \omega}{(\omega_3 + k)^3}- 
\frac{{|\psi|^2}}{(\omega_3 + k)^2}\Big)~\!\mu^2\!dz~\!,
$$
use (\ref{eq:J'dec1}) and a density argument. Next we put
$$
\lambda:= \inf_{\substack{\psi\in X \\ \psi\neq 0}}\frac{\displaystyle{\int_{\R^2} J'_0(\w)\psi\cdot\psi~\!dz}}
{\displaystyle{\int_{\R^2}(\omega_3 + k)^{-2}{{|\psi|^2}}~\mu^2\!dz}}\, ,
$$
and notice that $\lambda\ge 0$ by Lemma \ref{L:semidefpos}. On the other hand,
$\lambda$ is
achieved by 
Rellich theorem. Thus $\lambda>0$, because of Lemma \ref{L:tangsol}.
It follows that the energy functional $I: X \to \R$,
$$
I(\psi) =\frac{1}{2}\irn J'_0(\w)\psi\cdot \psi~\!dz - \int\limits_{\R^2} v \cdot \psi~\mu^2\!dz\, ,
$$
is weakly lower semicontinuous and coercive. Thus 
its infimum is achieved by a function $\f\in X$ which satisfies
\begin{equation}
\label{eq:quasi1}
\int\limits_{\R^2}J'_0(\w)\f\cdot \psi ~\!dz= \int\limits_{\R^2}  v \cdot \psi ~\mu^2\!dz\, , \quad \forall \ \! \psi \in X\, .
\end{equation}
If $\psi \in H^1$ we write
\[
\psi = (\P\psi^\top + \P\psi^\perp) + \eta \omega\, ,
\]
where $\eta=\psi\cdot\omega$, $\P\psi^\top \in \tange=\ker J'_0(U)$ is the orthogonal projection of $\P\psi=\psi-\eta\omega$ onto $\tange$
in the scalar product $(\cdot,\cdot)_*$ and
$\P\psi^\perp:= \psi -\P\psi^\top- \eta \omega \in X$. We use  \eqref{eq:SA} and \eqref{eq:J'dec2} to compute
\begin{gather*}
\int\limits_{\R^2}J'_0(\w)\f\cdot \P\psi^\top ~\!dz= \int\limits_{\R^2}J'_0(\w)\P\psi^\top\cdot \psi ~\!dz=0~\!,\\
\int\limits_{\R^2}J'_0(\w)\f\cdot (\eta\omega) ~\!dz= 
\irn \frac{\nabla(\f\cdot\omega)\cdot\nabla\eta}{(\omega_3 + k)^2} ~\!dz- 2k\irn
\frac{(\f\cdot\omega)\eta }{(\omega_3 + k)^3}~\mu^2\!dz =0\, ,
\end{gather*}
because $\f\cdot\omega\equiv 0$. Therefore,  (\ref{eq:quasi1}) gives
$$
\int\limits_{\R^2}J'_0(\w)\f\cdot \psi ~\!dz =\int\limits_{\R^2}J'_0(\w)\f\cdot \P\psi^\perp ~\!dz=  
\irn v \cdot \P\psi^\perp ~\!\mu^2dz= \irn v \cdot \psi ~\!\mu^2dz\,,
$$
as $v$ is orthogonal to $\tange \ni \P\psi^\top$ and to $\eta\omega$ in $L^2(\ovR^2,\R^3)$. We showed that 
$\f$ solves \eqref{eq:inverter0bis}, and thus the proof is complete.\QED

\begin{Lemma}\label{L:omeganonhom}
Let $f \in H^1(\ovR^2)$ be such that  $f \omega \in \tange^\perp$. There exists 
$\eta \in H^1(\ovR^2)$ such that 
$\eta\omega\in H^1\cap \tange^\perp_*\cap N_*$ and
\begin{equation}
\label{eq:inverter00}
J'_0(\w)(\eta\omega)=f\omega~\!\mu^2\quad\text{on $\R^2$.}
\end{equation}
\end{Lemma}

\proof
We introduce the space
\[
Y:= \Big\{\,\eta\in H^1(\ovR^2)~\Big|~ 
\int\limits_{\R^2}\frac{\eta}{(\omega_3 + {k})^3} ~\mu^2\!dz = \int\limits_{\R^2} \frac{\eta(\tau\cdot \omega)}{(\omega_3 + k)^3}~\mu^2\!dz\,  = 0\,,\  \forall \, \tau \in \tange \,  \Big\}\, ,
\]
so that $\eta\omega\in H^1\cap \tange^\perp_*\cap N_*$ for any $\eta\in Y$,
and the  energy functional $I:Y\to\R$,
\[
\begin{aligned}
I(\f) =&~ \frac{1}{2}\int\limits_{\R^2}J'_0(\w)(\eta \omega) \cdot (\eta \omega)~\!dz- \int\limits_{\R^2}  f\eta~\mu^2\!dz \\
=&~\frac12 \int\limits_{\R^2} \frac{|\nabla \eta|^2}{(\omega_3 + k)^2}~\!dz - k 
\int\limits_{\R^2}~\!\frac{|\eta|^2}{(\omega_3 + k)^3} \, \mu^2\!dz
- \int\limits_{\R^2} \eta f\, \mu^2\!dz\, ,
\end{aligned}
\]
compare with \eqref{eq:J'dec2}.
The functional $I$  is weakly lower semicontinuous with respect to the $H^1(\ovR^2)$ topology and 
coercive by Remark \ref{R:spectrum}. Thus 
its infimum is achieved by a function $\eta\in Y$. To conclude, argue as in the proof of Lemma \ref{L:tangnonhom}  to
show that $\eta$ solves (\ref{eq:inverter00}).
\QED

\paragraph{Proof of Lemma \ref{L:Fredholm}.} 
Since $J'_0(\w_{\!q})=q_3^{-2}J'_0(\w)$, we can assume that $q=e_3$, that is, $\w_{\!q}=\w$.
We take any $v\in \tange^\perp$, and write
$$
v = \P v + (v\cdot\omega)\omega\, ,
$$ 
where $\P v= v-(v\cdot\omega)\omega$, as before. Since
$\P v\in \tange^\perp$, by Lemma \ref{L:tangnonhom} there exists a unique $
\hat\f\in  H^1\cap \tange^\perp_*$ such that
$\hat\f\cdot\omega\equiv 0$ on $\R^2$ and
$$
\irn J'_0(\w)\hat\f\cdot \psi~\!dz\, =\irn \P v\cdot\psi~\!\mu^2\!dz\, , \quad\text{for any $\psi\in H^1$.}
$$
Next, notice that $(v\cdot\omega)\omega\in \tange^\perp$, so we can use Lemma 
\ref{L:omeganonhom} to find $\eta \in H^1(\ovR^2)$ such that 
$\eta\omega\in H^1\cap \tange^\perp_*\cap N_*$ solves
$$
\irn J'_0(\w)(\eta\omega)\cdot \psi~\!dz=
\irn (v\cdot\omega)(\psi\cdot\omega)~\!\mu^2\!dz, \quad\text{for any $\psi\in H^1$.}
$$
The function $\f_v=\hat\f+\eta\omega$ solves (\ref{eq:inverter0}).

\medskip

To conclude the proof we have to show that if $v\in {\cal C}^{m}(\ovR^2, \R^3)$ then $\f_v \in {\cal C}^{2+m}(\ovR^2, \R^3)$. 
Since $\omega \in {\cal C}^{\infty}(\ovR^2,\R^3)$ and $\omega_3 + k$ is bounded and bounded away from zero, 
$\f_v$ solves a linear system of the form
 $$
 -\Delta \f_v=A(z)\f_v+ B(z)\nabla \f_v+\mu^2(\omega_3+k)^2 v\, ,
 $$
 for certain smooth matrices on $\ovR^2$. A standard bootstrap argument and Schauder regularity theory
 plainly imply that $\f_v \in {\cal C}^{2+m}_{loc}(\R^2, \R^3)$. 
 The function $z\mapsto \f_v(z^{-1})$ satisfies a linear system of the same kind, hence $\f_v \in {\cal C}^{2+m}(\ovR^2, \R^3)$, as desired.  \QED

\section{The perturbed problem}
\label{S:fdr}

In this Section we perform the finite dimensional reduction and prove Theorems \ref{T:existence}, \ref{T:main2}.
 By the results in the Appendix, any 
critical point of the ${\cal C}^2$-functional $E_\eps:{\cal C}^{2}(\ovR^2,\H^3)\to \R$, 
$$
E_\eps(u):=
\frac12\irn  u_3^{-2}|\nabla u|^2~\! dz- k\irn u_3^{-2}e_3\cdot \partial_x u\wedge\partial_yu~\!dz
+2\eps ~\!V_\phi(u)=E_0(u)+2\eps ~\!V_\phi(u)
$$
(notice that $E_0=E_{\text{hyp}}$, compare with (\ref{eq:ene})), solves
\begin{equation*}\tag{$\mathcal P_\eps$}
\label{eq:eps_problem}
\Delta u - 2 u_3^{-1} {{G}(\nabla u)} =2 u_3^{-1}(k+\eps\phi(u)) ~\!\partial_x u\wedge \partial_y u\quad \text{on $\R^2$}
\end{equation*}
and has mean curvature $(k+\eps \phi)$, apart from a finite set of branch points. 

Due to the action of the M\"obius transformations and of the hyperbolic translations, for any $u\in {\cal C}^{2}(\ovR^2,\R^3)$ we have the identities
\begin{gather}\label{eq:Econfinv}
E_\eps'(u)(z^h\nabla u)=0~,~~ E_\eps'(u)(iz^h\nabla u)=0~,~~\text{for $h=0,1,2$\,,\ \ $\eps\in\R$\,,}\\
\label{eq:Etrinv}
E'_0(u) e_1 = 0\,, \quad E'_0(u)e_2=0\, ,\quad E_0'(u) u = 0~\!.
\end{gather}
Now we prove that
\begin{equation}
\label{eq:oggi2}
E_\eps(\w_{\!q})= E_0(\w)- 2\eps F^{\phi}_k(q)\, ,
\end{equation}
where $F^{\phi}_k$ is the Melnikov type function in \eqref{eq:Mel}. The above mentioned invariances give $E_0(\w_{\!q})=E_0(\w)$. 
Since the hyperbolic ball $B^\H_{\rho_k}(q)$ coincides with the Euclidean ball of radius $q_3r_{\!k}$
about the point $q^k:=(q_1,q_2,kr_{\!k}q_3)$, the divergence theorem gives
$$
F^{\phi}_k(q)=\int\limits_{B^\H_{\!\rho_{k}}(q)}\phi(p)~\!d\H^3_p\, =\! \int\limits_{B_{q_3r_{\!k}}(q^k)}p_3^{-3}\phi(p)~\!dp\,=\!
\int\limits_{\partial B_{q_3r_{\!k}}(q^k)}Q_\phi(p)\cdot \nu_p\, .
$$
Here $Q_\phi\in {\cal C}^1(\R^3_+,\R^3)$ is any vectorfield such that $\text{div}Q_\phi(p)=p_3^{-3}\phi(p)$
and $\nu_p$ is the outer normal to $\partial B_{q_3r_{\!k}}(q^k)$ at $p$. The function $\w_{\!q}$ in (\ref{eq:wq})
parameterizes the Euclidean sphere $\partial B{q_3r_{\!k}(q^k)}$. Since $\partial_x \w_{\!q}\wedge \partial_y\w_{\!q}$
is inward-pointing, we have
\begin{equation}
\label{eq:FV}
F^{\phi}_k(q)=-\int\limits_{\R^2}Q_\phi(p)\cdot \partial_x \w_{\!q}\wedge \partial_y\w_{\!q}~\!dz=-V_\phi(\w_{\!q}),
\end{equation}
and (\ref{eq:oggi2}) is proved. Before going further, let us  show that the existence of critical points for $F^{\phi}_k$ is a necessary condition 
for the conclusion in Theorem \ref{T:existence}. 

\begin{Theorem}
\label{T:side1}
Let $k>1$, $\phi\in {\cal C}^{1}(\H^3)$. Assume that 
there exist  sequences $\eps_h\subset\R\!\setminus\!\{0\}$, $\eps_h\to 0$, $u^h\in {\cal C}^2(\ovR^2,\H^3)$
and a point $q\in \H^3$ such that 
$u_h$  solves $(\mathcal P_{\eps_h})$, and $u^h\to \w_{\!q}$ in ${\cal C}^1(\ovR^2,\H^3)$.
Then $q$ is a stationary point for $F^{\phi}_k$.
\end{Theorem}

\proof
The function $u^h$ is a stationary point for the energy functional 
$E_{\eps_h}=E_0+2\eps_h V_\phi$.  From (\ref{eq:Etrinv}) we have
$V'_\phi(u^h)e_j=0$ for $j=1,2$ and $V'_\phi(u^h)u^h=0$. We can plainly pass to the limit to obtain
$V'_\phi(\w_{\!q})e_j=0$  for $j=1,2$ and $V'_\phi(\w_{\!q})\w_{\!q}=0$. To conclude, use (\ref{eq:FV}) and recall that 
$\partial_{q_j} \w_{\!q} =e_j$  for $j=1,2$, and $\partial_{q_3} \w_{\!q}=\w=q_3^{-1}(\w_{\!q}-q_1e_1-q_2e_2)$. 
\QED

\medskip
Now we fix $m\in(0,1)$. The  operator $J_\e:{\cal C}^{2+m}(\ovR^2,\H^3)\to {\cal C}^m(\ovR^2,\R^3)$ defined by 
$$
J_\e(u)=-\div(u_3^{-2}\nabla u)-u_3^{-3}|\nabla u|^2e_3+2(k+ \e\phi)u_3^{-3}{\partial_{x\!}u}\wedge {\partial_{y\!}u}\,,
$$
is related to the differential of $E_\eps$ via the identity
\begin{equation}
\label{eq:E'J}
E_\e'(u)\f = \irn J_\e(u)\cdot \f~\!dz~,\quad u\in {\cal C}^m(\ovR^2,\H^3)\,,\ \f \in {\cal C}^m(\ovR^2,\R^3)\,.
\end{equation}

\begin{Remark}
\label{R:aggiunta}
Since $E_\eps$ is of class ${\cal C}^2$ and
$$
E''_\e(u)[\f,\psi]=\displaystyle\int_{\R^2} J'_\e(u)\psi\cdot\f~\!dz\,,
$$
then the quadratic form in the right hand side is a self-adjoint form on $H^1$.
\end{Remark}

We are in position to state and proof the next lemma, which is the main step towards the proofs of Theorems \ref{T:existence}, \ref{T:main2}. 

\begin{Lemma}[Dimension reduction]
\label{L:findimred}
Let $\Omega \Subset \H^3$ be an open set. There exists $\hat \e >0$ and a unique ${\cal C}^1$-map
\[
[-\hat \e, \hat \e]\times \overline \Omega \to {\cal C}^{2+m}(\ovR^2, \H^3)\,, \quad (\e, q) \mapsto u^\e_q\,,
\] 
such that the following facts hold:
\begin{enumerate}
\item[$i)$] $u^{\e}_q$ parameterizes an embedded $\S^2$-type surface, and  $u^{0}_q = \w_{\!q}$\,;
\item[$ii)$] $u^\e_q - \w_{\!q} \in \tange^\perp\cap \mathcal{C}^{2+m}(\ovR^2, \R^3)$ and $E'_\e(u^\e_q)\f=0$ for any $\f\in  \tange^\perp\cap \mathcal{C}^{0}(\ovR^2, \R^3)$\,;
\item[$iii)$] for any $\eps\in [-\hat \e, \hat \e]$, the manifold $\{\,u^\eps_q~|~ q\in \Omega\,\}$ is a natural constraint for $E_\eps$, that is,
if $\nabla_{\!\!q}  E_\e(u^\e_{q^\eps}) = 0$ for some $q^\e\in\Omega$, then $u^{\e}_{q^\e}$ is a $(k + \e\phi)$-bubble\,;
\item[$iv)$] 
$\|E_\e(u^{\e}_q) - E_\e(\w_{\!q})\|_{{\cal C}^1\left(\overline\Omega\right)} = o (\e)$ as $\e \to 0$, uniformly on $\overline\Omega$\,.
\end{enumerate}
\end{Lemma}

\proof
To shorten the notation, we put $\mathcal{C}^m:=\mathcal{C}^m(\ovR^2, \R^3)$. For $s \geq0$ and $\delta >0$ we write
\[
\begin{aligned}
\Omega_s := \{\, p \in \H^3 \ |\ \text{dist}(p, \Omega) < s\, \}\, ,\, \ \text{and}\ \ \mathcal{U}_\delta := \{\, \nu \in {\cal C}^{2+m} \ |\ |\nu(z)|<\delta \text{ for every }z \in \R^2\, \}\, .
\end{aligned}
\]
We fix $s>0$, $\delta=\delta(s)>0$ such that $\overline \Omega_{2s} \subset \H^3$ and  $(\w_{\!q}+ \nu)\cdot e_3 >0$ for  $q \in \Omega_{2s}$, $\nu \in \mathcal{U}_\delta$. 

\newcommand\T{\mathbb T}
We define 
\begin{equation}\label{eq:taugammadef}
\begin{array}{lll}
\tau_1:=c_0\partial_x\omega~\!,&\tau_3:= c_0\sqrt2 z\nabla \omega~\!,&\tau_5:= c_0z^2\nabla \omega~\!,\\
\tau_2:=c_0\partial_y\omega~\!,&\tau_4:=c_0 \sqrt2  iz\nabla \omega~\!,&\tau_6:=c_0iz^2\nabla \omega~\!,
\end{array}
\qquad \gamma:=2c_0(k\omega+e_3)\, ,
\end{equation}
where $c_0:=\sqrt{\frac{3}{2^4\pi}}$ is a normalization constant. Thanks to \eqref{eq:tangbasis},
(\ref{eq:TUZ}), we have 
$$
\tange=\langle \tau_1,\dots\tau_6\rangle\oplus ~\{\, (\alpha\cdot\gamma)~\!\omega~|~\alpha\in\R^3\, \}~\!.
$$
Trivially, $\tau_j\cdot\omega\equiv 0$ on $\R^2$. Elementary computations give 
$$
\int\limits_{\R^2}{\tau_i \cdot \tau_{j}}~\mu^2\!dz  = \delta_{ij}~\!,\quad 
 \int\limits_{\R^2}{\gamma_h \gamma_\ell} ~\mu^2\!dz = 0\quad \text{if $h\neq\ell$}~\!,
  $$
 for $i,j \in\{1,\dots,6\}$, $h,\ell\in\{1,2,3\}$, and moreover
$$
 \irn \gamma_1^2~\mu^2\!dz=\irn \gamma_2^2~\mu^2\!dz= k^2~,\quad 
 \irn \gamma_3^2~\mu^2\!dz=k^2+3~\!.
  $$
\paragraph{Construction of $\bf{u^\eps_q}$ satisfying \bf{{\em i})}, \bf{{\em ii})}.} By our choices of $s$ and $\delta$, the functions
$$
\begin{aligned}
&\mathcal{F}_1(\e,q; \nu, \xi,\alpha) := \mu^{-2}J_\e(\w_{\!q} + \nu)- 
\sum_{j=1}^6 \xi_j\tau_j- (\alpha\cdot\gamma) ~\!\omega~\in {\cal C}^m \,, \\
&\mathcal{F}_2(\e, q; \nu, \xi,\alpha ) := \Big(\ 
\int\limits_{\R^2}{\nu \cdot \tau_1}\ \mu^2\!d z\ ,\dots, \int\limits_{\R^2}{\nu \cdot \tau_6}\ \mu^2\!d z~;\ \
\int\limits_{\R^2}\gamma~\!(\nu\cdot \omega)\ \mu^2\!d z\ \Big)~\in \R^6 \times \R^3 \,, 
\end{aligned}
$$
are well defined and continuously differentiable on $\R \times \Omega_{2s}\times \mathcal{U}_\delta \times (\R^6 \times \R^3)$. Thus
$$
\mathcal{F} := (\mathcal{F}_1, \mathcal{F}_2) : \R \times \Omega_{2s}\times \mathcal{U}_\delta \times (\R^6 \times \R^3) \to  {\cal C}^m\times (\R^6 \times \R^3)
$$
is of class ${\cal C}^1$ on its domain. Notice that $\F(0, q; 0, 0, 0) = 0$ for every $q \in \Omega_{2s}$ because $J_0(\w_{\!q}) = 0$. Now we solve the equation $\F(\e, q; \nu, \xi,\alpha) =0$ in a neighborhood of $(0, q; 0, 0, 0)$ via the implicit function theorem. 
Let
\[
\L := (\L_1, \L_2) : {\cal C}^{2+m} \times (\R^6 \times \R^3) \to {\cal C}^m \times (\R^6 \times \R^3)
\]
given by
\[
\begin{aligned}
&\mathcal{L}_1(\f; \zeta, \beta) := \mu^{-2}J'_0(\w_{\!q})\f - \sum_{j=1}^6 \zeta_j \tau_j - (\beta\cdot \gamma)~\!\omega\, ,\\
&\mathcal{L}_2(\f; \zeta, \beta ) := \mathcal{L}_2(\f)=
\Big(\ 
\int\limits_{\R^2}\f \cdot \tau_1\ \mu^2\!d z\ ,\dots, \int\limits_{\R^2}\f \cdot \tau_6\ \mu^2\!d z\ 
~;~
\int\limits_{\R^2}\gamma~\!(\f \cdot \omega)\ \mu^2\!d z\ \Big)\, ,
\end{aligned} 
\] 
so that 
$\L=(\L_1,\L_2)$
 is  the differential of $\F(0,q; \cdot , \cdot, \cdot)$ evaluated in $(\nu, \xi,\alpha) = (0,0,0)$. 

To prove that $\mathcal{L}$ is injective we assume that $\L(\f, \zeta, \beta) = 0$ and put
$$v= \mu^{-2}J'_0(\w_{\!q})\f \in \tange\,.$$ 
From \eqref{eq:SA} we find
\[
\int\limits_{\R^2}|v|^2\ \mu^2\!d z = \int\limits_{\R^2} \big(\mu^{-2}J'_0(\w_{\!q})\f)\cdot v\ \mu^2\! dz  = 
\int\limits_{\R^2} J'_0(\w_{\!q})\f\cdot v \ dz=\int\limits_{\R^2} J'_0(\w_{\!q})v\cdot \f \ dz  = 0\, ,
\]
which implies $J'_0(\w_{\!q})\f=0$, that is, $\f \in \tange$.  On the other hand, $\f \in T_{\!\w\!}Z^\perp$
because $\L_2(\f) = 0$. Thus $\f = 0$ and therefore also $\beta = \zeta = 0$.

To prove that $\L$ is surjective fix  $v \in {\cal C}^m$ and $(\theta, b) \in \R^6 \times \R^3$. We have to find  $\f \in {\cal C}^{2+m}$ and $(\zeta, \beta) \in \R^6 \times \R^3$ such that
$\L_1(\f;\zeta,\beta)=v$ and $\L_2(\f)=(\theta,b)$. To this goal we introduce the
minimal distance projection
$$
P^\top: L^2(\ovR^2,\R^3)\to \tange\, , \quad w \mapsto P^{\!\top}\!w\, ,
$$
so that $\L_2(w)$ is uniquely determined by $P^\top \!w$, and vice-versa.
We
find $\zeta_j$ and 
$\beta$ so that
$$
\sum_{j=1}^6 \zeta_j \tau_j + (\beta\cdot\gamma)\omega=-P^{\!\top}\!v\, .
$$
Then, we use Lemma \ref{L:Fredholm} to find $\widehat\f\in {\cal C}^{2+m}\cap \tange^\perp_*\cap N_*$ such that
$$
J'_0(\w_{\!q})\widehat\f=(v-P^{\!\top}\!v)~\!\mu^2.
$$
Finally, we take the unique tangent direction $\f^\top\!\in \tange$ such that 
$\L_2(\f^\top)=(\theta,b)-\L_2(\widehat\f)$.
The triple $(\f^\top\!+\widehat\f,\zeta,\beta)$ satisfies $\L(\f^\top\!+\widehat\f,\zeta,\beta)=(v;\theta,b)$ and surjectivity is proved.
We are in the position to apply the implicit function theorem to $\F$, for any fixed $q \in \Omega_{2s}$. In fact, thanks to a standard compactness argument, we get that there exist $ \e' >0$
and uniquely determined ${\cal C}^1$ functions 
\[
\begin{aligned}
&\nu:(-\eps',\eps') \times \Omega_s \to \mathcal{U}_\delta & \quad& \alpha: (-\eps',\eps') \times \Omega_s \to \R^3 & \quad&\xi: (-\eps',\eps') \times \Omega_s \to \R^6\\
&\nu: (\e, q) \mapsto \nu^{\e}_q & &\alpha: (\e, q) \mapsto \alpha^{\e}(q)  &&\xi: (\e, q) \mapsto \xi^{\e}(q)
\end{aligned}
\]
such that
\begin{equation}\label{eq:implsol}
\nu^{0}_q \equiv 0, \qquad \alpha^0(q) = 0\, , \qquad \xi^{0}(q) = 0\, ,  \qquad\mathcal{F}(\e, q; \nu^{\e}_q, \xi^{\e}(q), \alpha^{\e}(q)) = 0\, .
\end{equation}
By \eqref{eq:implsol}, the ${\cal C}^1$ function $(-\eps',\eps') \times \Omega_s \to {\cal C}^{2+m}(\ovR^2, \H^3)$,
\[
(\e, q) ~\mapsto ~u^\e_q := \w_{\!q}+ \nu^{\e}_q=
~\big(q_3\w+q_1e_1+q_2e_2\big)+\nu^{\e}_q\,,
\]
 satisfies $i)$, if $\e'$ is small enough. Further, using
  (\ref{eq:E'J}) (see also Lemma \ref{A:funct}) we rewrite the last identity in (\ref{eq:implsol}) as 
\begin{equation}\label{eq:implprop}
\begin{aligned}
&\begin{aligned}
E'_\eps(u^\e_q)\f=&\irn J'_\eps(\w_{\!q}+ \nu^{\e}_q)\cdot \f~\!dz\\
=& \sum_{j=1}^6 {\xi_j^{\e}(q)}\irn\tau_j\cdot\f~\mu^2\!dz + 
\irn (\alpha^\e(q)\cdot \gamma)(\omega\cdot\f)~\mu^2\!dz\quad \text{$\forall\,\f\in {\cal C}^0$,}
\end{aligned}
\\
&\int\limits_{\R^2}\nu^{\e}_q \cdot \tau_j~ \mu^2\!dz = 0\,,\quad\text{$\forall\, j\in\{\, 1, \ldots, 6\,\}$}\, , \quad 
\ \int\limits_{\R^2}\gamma_\ell(\nu^{\e}_q \cdot \omega)\ \mu^2\!dz= 0\,,\quad\text{$\forall\, \ell \in \{\, 1,2,3\,\}$}\, .
\end{aligned}
\end{equation}
In particular, claim $ii)$ holds true.

\paragraph{Proof of $\bf{iii)}$.}
As a straightforward consequence of (\ref{eq:implprop}) we have that 
\[
\int\limits_{\R^2}\partial_{q_i}\nu^{\e}_q \cdot \tau_j\ \mu^2\!dz = 0\,, \quad \quad 
\int\limits_{\R^2}\gamma_\ell~\!(\partial_{q_i}\nu^{\e}_q\cdot \omega)\ \mu^2\!dz= 0\,,
\]
hence $E'_\e(u^\e_q)\partial_{q_i}\nu^{\e}_q = 0$ for any $ i=1,2,3$. We infer the identities
\begin{equation}\label{eq:fE'rel}
\begin{aligned}
&\partial_{q_i}E_\e(u^\e_q) = E'_\e(u^\e_q)(e_i + \partial_{q_i}\nu^{\e}_q) = E'_\e(u^\e_q)e_i\, , \quad i=1,2\, ,\\
&\partial_{q_3}E_\e(u^\e_q) = E'_\e(u^\e_q)(\w + \partial_{q_3}\nu^{\e}_q) = E'_\e(u^\e_q)\w\, .
\end{aligned}
\end{equation}
Now, from  \eqref{eq:omega3}, \eqref{eq:taugammadef} and \eqref{eq:implprop} we find
\begin{gather*}
2c_0e_1 = \tau_1 - \tau_5 + k^{-1}\gamma_1 \omega\, , \quad 2c_0e_2 = \tau_2 + \tau_6 + k^{-1}\gamma_2 \omega\, ,\quad
2c_0\w =k r_k(\sqrt{2}\tau_3 +k^{-1}\gamma_3 \omega)\, ,\\
E'_\e(u^\e_q)\tau_j  = \xi^{\e}_j(q)\, ,\qquad
E'_\e(u^\e_q)(\gamma_\ell \omega) = (k^2+3\delta_{\ell3})\alpha_\ell^{\e}(q)\,,
\end{gather*}
for any $j =1, \ldots, 6$, $\ell=1,2,3$.
Thus by \eqref{eq:fE'rel} we get  
\begin{equation}
\label{eq:fgrad}
2c_0\nabla_{\!q}E_\eps(u^\eps_q)=M_k\xi^\eps(q)+\Theta_k\alpha^\eps(q)~\!,
\end{equation}
where $M_k$ and $\Theta_k$ are constant matrixes, namely
$$
M_k=\left(\begin{array}{cccccc}
1&0&0&0&-1&0\\
0&1&0&0&0&-1\\
0&0&{\sqrt{2}kr_{k}}&0&0&0
\end{array}\right)~,\qquad 
\Theta_k=\left(\begin{array}{ccc}
k&0&0\\
0&k&0\\
0&0&(k^2+3)r_k
\end{array}\right)~.
$$
On the other hand, from \eqref{eq:Econfinv} and using $\nabla \w_{\!q} = r_kq_3 \nabla \omega$ we obtain
\begin{equation}
\label{eq:nohorror}
-q_3r_k~\!\xi^{\e}_j(q)=E_\e'(u^{\e}_q)(\tau^\eps_j(q))\,,
\end{equation}
where, in the spirit of (\ref{eq:taugammadef}), we have putted
$$
\begin{array}{lll}
\tau^\eps_1(q):=c_0\partial_x\nu^\eps_q~\!,&\tau^\eps_3(q):=c_0\sqrt{2}z\nabla \nu^\eps_q~\!,&\tau^\eps_5(q):= c_0z^2\nabla \nu^\eps_q~\!,\\
\tau^\eps_2(q):=c_0\partial_y\nu^\eps_q~\!,&\tau^\eps_4(q):= c_0\sqrt2  iz\nabla \nu^\eps_q~\!,&\tau^\eps_6(q):=c_0iz^2\nabla \nu^\eps_q~\!.
\end{array}
$$
Notice that 
\begin{equation}
\label{eq:Holder}
\irn|\tau^\eps_j(q)|^2~\mu^2\!dz\le 2\irn |\nabla_{\!z}\nu^\eps_q|^2~\mu~\!dz\le 2~\|\nu^\eps_q\|_{{\cal C}^1}^2\irn\mu^3~\!dz=o(1),
\end{equation}
as $\eps\to 0$, uniformly on $\overline \Omega$, see (\ref{eq:Cm_norm}).

For sake of clarity, we make now some explicit computation. We denote by $\sigma_{\ell h}$ the entries of the $3\times 6$ constant matrix $\Theta^{-1}_kM_k$,
and introduce the $6\times 6$ matrix $A^\eps(q)=(a^\eps_{jh}(q))_{j,h=1,\dots,6}$, whose entries
are given by
$$
a^\eps_{jh}(q)=\irn \tau_h\cdot \tau^\eps_j(q)~\mu^2\!dz-
\sum_{\ell=1}^3 \sigma_{\ell h}  \irn \gamma_\ell (\omega\cdot \tau^\eps_j(q))~\mu^2\!dz~\!.
$$

Since $\tau^\eps_j\mu\to 0$ in $L^2(\R^2,\R^3)$ by (\ref{eq:Holder}), then $A^\eps\to 0$ uniformly on compact subsets of $(-\eps',\eps')\times\Omega_s$. In particular,
if $\hat\eps\in(0,\eps')$ is small enough, then the determinant of the $6\times 6$ matrix $(A^\eps(q)+q_3r_k\text{Id})$ is uniformly bounded away from $0$
on $[-\hat\eps,\hat\eps]\times\overline{\Omega}$.

\medskip

Assume that $\nabla_{\!q}E_\eps(u^\eps_{q^\eps})=0$ for some $\eps\in [-\hat\eps,\hat\eps]$, $q^\eps\in \Omega$.
From (\ref{eq:fgrad}) we obtain   $\alpha^\eps(q^\eps)=-\Theta_k^{-1}M_k\xi^\eps(q^\eps)$. Thus 
\eqref{eq:implprop} and \eqref{eq:nohorror} give
$$
-q^\e_3 r_k~\!\xi^\eps(q^\eps)= A^\eps(q^\e)\xi^\eps(q^\eps),
$$
and hence $\xi^\eps(q^\eps)=0$, because the matrix $(A^\eps(q^\e)+q^\e_3r_k\text{Id})$ is invertible. But then  (\ref{eq:fgrad}) and $\nabla_{\!q}E_\eps(u^\eps_{q^\eps})=0$
imply that $\alpha^\eps(q^\eps)=0$ as well, hence $E'(u^\eps_{q^\eps})=0$ by (\ref{eq:implprop}).

\paragraph{Proof of $\bf{iv)}$.} 
The function $(\e, q) \mapsto \nu^{\e}_q$ is of class ${\cal C}^1$, and in particular  $\partial_\eps\nu^{\e}_q$ is uniformly bounded in ${\cal C}^2$ 
for $(\eps,q)\in[-\hat\eps,\hat\eps]\times\overline\Omega$. Thus Taylor expansion formula for
$$
\eps \mapsto E_\eps(u^\eps_q)-E_\eps(\w_{\!q})=E_0(u^\eps_q)-E_0(\w_{\!q})+2\eps\big(V_\phi(u^\eps_q)-V_\phi(\w_{\!q}))
$$
gives $E_\eps(u^\eps_q)-E_\eps(\w_{\!q})=o(\eps)$ as $\eps\to 0$, uniformly on $\overline\Omega$.

Now we estimate $\nabla_{\!q}(E_\eps(u^\eps_q)-E_\eps(\w_{\!q}))$. We use (\ref{eq:Etrinv}), \eqref{eq:fE'rel} to obtain, for $j=1,2$, 
$$
\begin{aligned}
\partial_{q_j}(E_\eps(u^\eps_q)-E_\eps(\w_{\!q}))=&~\big(E'_0(u^\eps_q)e_j-E'_0(\w_{\!q})e_j\big)+2\eps\big(V'_\phi(u^\eps_q)e_j-V'_\phi(\w_{\!q})e_j\big)\\
=&~2\eps\big(V'_\phi(u^\eps_q)e_j-V'_\phi(\w_{\!q})e_j\big)=o(\eps)\,,
\end{aligned}
$$
because $\|u^\eps_q-\w_{\!q}\|_{{\cal C}^{2+m}}=o(1)$ and $V_\phi$ is a ${\cal C}^1$-functional.  

To handle the derivative with respect to $q_3$ we first argue as before to get
$$
\begin{aligned}
\partial_{q_3}(E_\eps(u^\eps_q)-E_\eps(\w_{\!q}))=&~\big(E'_0(u^\eps_q)\w-E'_0(\w_{\!q})\w\big)+2\eps\big(V'_\phi(u^\eps_q)\w-V'_\phi(\w_{\!q})\w\big)\\
=&~E'_0(u^\eps_q)\w+ o(\eps)~\!,
\end{aligned}
$$
uniformly on $\overline\Omega$.
Next, from $q_3\w=u^\eps_q-(q_1e_1+q_2e_2)-\nu^\eps_q$ and (\ref{eq:Etrinv}) we obtain
\[
\begin{aligned}
q_3E'_0(u^\eps_q)\w=&E'_0(u^\eps_q)(u^\eps_q-(q_1e_1+q_2e_2)-\nu^\eps_q)\\
=&-E'_0(u^\eps_q)\nu^\eps_q
=-E'_\eps(u^\eps_q)\nu^\eps_q+2\eps V'_\phi(u^\eps_q)\nu^\eps_q=2\eps V'_\phi(u^\eps_q)\nu^\eps_q
\end{aligned}
\]
because of \eqref{eq:implprop}. Since $\nu^\eps_q\to 0$ in ${\cal C}^{2+m}$ we infer that
$E'_0(u^\eps_q)u^\eps_q=o(\eps)$ uniformly on $\overline\Omega$ as $\eps\to 0$, which concludes the proof.
\QED

\paragraph{Proof of Theorem \ref{T:existence}.} 
Take an open set $\Omega\Subset\R^3_+$ containing the closure of $A$,  let $u^\eps_q$ be the function given by 
Lemma \ref{L:findimred} and notice that, by \eqref{eq:FV},
$E_\eps(\w_{\!q})=E_0(\w_{\!q})- 2\e F^\phi_k(q)$. Thus for $\eps\in [-\hat\eps,\hat\eps], \eps\neq 0$ we can estimate
$$
\Big\|\frac{1}{2\eps}\big(E_\eps(u^\eps_q)-E_0(\w_{\!q})\big)+F^\phi_k(q)\Big\|_{{\cal C}^1(\overline A)}=\frac{1}{2|\eps|}\|E_\eps(u^\eps_q)-E_\eps(\w_{\!q})\|_{{\cal C}^1(\overline A)}=o(1)\,,
$$
uniformly on $\overline \Omega$ by $iv)$ in Lemma \ref{L:findimred}. Recalling the definition of stable critical point presented in Subsection \ref{SS:stable}, we infer that for any $\eps\approx 0$ the function $\frac{1}{2\eps}\big(E_\eps(u^\eps_q)-E_0(\w_{\!q})\big)$ has a critical point $q^\e \in A$,
to which corresponds the embedded $(k+\e\phi)$-bubble $u^\eps:=u^\eps_{q^\eps}$ by $iii)$ in Lemma \ref{L:findimred}.
The continuity of $(\eps,q)\mapsto u^\eps_q$ gives the continuity of $\eps\mapsto u^\eps$. 

The last conclusion in Theorem \ref{T:existence} follows via a simple compactness argument
and thanks to Theorem \ref{T:side1}. \QED

\paragraph{Proof of Theorem \ref{T:main2}.} 
Recalling that $q^k:=(q_1,q_2,kr_{\!k}q_3)$, we write 
$$
F^{\phi}_k(q)=\int\limits_{B_{r_{\!k}}(0)}(p_3+kr_{\!k})^{-3}\phi(q_3p+q^k)~\!dp~\!.
$$
Since $r_k \to 0$ and $kr_k=k(k^2-1)^{-1/2}\to 1$ as $k\to \infty$, we infer that $q^k\to q$ uniformly on compact sets of $\R^3_+$ and
$$
\frac{3}{4\pi r_{\!k}^3} F^{\phi}_k \to \phi \qquad\text{as $k\to\infty$}\,,
$$
uniformly on $\overline\Omega$. Next, we easily compute
$$
\begin{gathered}
\partial_{q_j} F^{\phi}_k(q)=\int\limits_{B_{r_{\!k}}(0)}(p_3+kr_{\!k})^{-3}\partial_{q_j}\phi(q_3p+q^k)~\!dp\,,\qquad j=1,2\,,\\
\partial_{q_3} F^{\phi}_k(q)=\int\limits_{B_{r_{\!k}}(0)}(p_3+kr_{\!k})^{-3}\nabla\phi(q_3p+q^k)\cdot(p+kr_{\!k}e_3)~\!dp\,,
\end{gathered}
$$
and thus we obtain, by the same argument,
$$
\frac{3}{4\pi r_{\!k}^3} \nabla F^{\phi}_k \to \nabla \phi\qquad\text{as $k\to\infty$,}
$$
uniformly on $\overline\Omega$. It follows that for $k$ large enough,
$F^{\phi}_k$ has a stable critical point in $\Omega \Subset \H^3$,
since having a stable critical point is a ${\cal C}^1$-open condition. Thus
Theorem \ref{T:main} applies and gives the conclusion of the proof. \QED

\appendix
\section{\!\!\!\!\!\!ppendix}
\label{A:appendix}

Let $K\in {\cal C}^0(\H^3)$. Take any vectorfield $Q_K\in{\cal C}^1(\R^3_+, \R^3)$ such that
$\div Q_K(p)=p_3^{-3}K(p)$ for any $p\in\R^3_+$ (here $\div=\sum_j\partial_j$ is the Euclidean divergence). The
functional
$$
V_K(u):=\irn Q_K(u)\cdot\partial_x u\wedge\partial_yu~\!dz\,,\qquad u\in{\cal C}^1(\ovR^2,\H^3)\,,
$$
measures the signed (hyperbolic) volume enclosed by the surface $u$, with respect to the weight $K$. In fact,  
if $u$ parameterizes the boundary of a smooth open set $\Omega\Subset \R^3_+$ and if $\partial_xu\wedge\partial_y u$
is inward-pointing, then the divergence theorem gives
$$
V_K(u)=-\int\limits_{\partial \Omega} Q_K(u)\cdot \nu~\!du=-\int\limits_\Omega p_3^{-3}K~\!dp=-\int\limits_\Omega K~\!d \H^3\, .
$$
Clearly, the functional $V_K$ does not depend on the choice of the vectorfield $Q$. Notice that if $K\equiv k$ is constant, then
$$
V_k(u)=-\frac{k}{2} \irn u_3^{-2}e_3\cdot \partial_x u\wedge\partial_yu~\!dz~,\quad u\in {\cal C}^1(\ovR^2,\H^3)~\!.
$$
In the next Lemma we collect few simple remarks about the energy functional 
\begin{equation}
\label{eq:energy_gen}
{E}(u)=\frac12\irn  u_3^{-2}|\nabla u|^2~\! dz+2V_{\!K\!}(u)\, .
\end{equation}

\begin{Lemma}\label{A:funct}
\label{L:volume}
Let $K\in {\cal C}^0(\H^3)$. 
\begin{itemize}
\item[$i)$] The functional $E:{\cal C}^1(\ovR^2,\H^3)\to\R$ is of class ${\cal C}^1$, and its differential is given by
$$
E'(u)\f=\irn(u_3^{-2}\nabla u\cdot\nabla\f-u_3^{-3}|\nabla u|^2e_3\cdot\f)~\!dz+2\irn u_3^{-3}K(u)\f\cdot \partial_x u\wedge\partial_yu~\!dz~;
$$
\item[$ii)$] If $u\in {\cal C}^2(\ovR^2,\H^3)$, then $E'(u)$ extends to a continuous form on 
${\cal C}^0(\ovR^2,\R^3)$\,, namely
$$
E'(u)\f=\irn(-\div(u_3^{-2}\nabla u)-u_3^{-3}|\nabla u|^2e_3+2u_3^{-3}K(u)\partial_x u\wedge\partial_yu)
\cdot\f~\!dz~;
$$
\item[$iii)$] If $K\in  {\cal C}^1(\H^3)$, then $E$ is of class ${\cal C}^2$ on ${\cal C}^1(\ovR^2,\H^3)$.
\end{itemize}
\end{Lemma}

In the next Lemma we show that critical points for ${E}$ are in fact hyperbolic $K$-bubbles.

\begin{Lemma}
\label{L:conformal}
Let $K\in {\cal C}^0(\H^3)$ and  let $u\in {\cal C}^2(\ovR^2,\H^3)$ be a nonconstant critical point for ${E}$. Then $u$
is conformal, that is,
$$
|\partial_x u|=|\partial_y u|\, ,~~ \partial_x u\cdot \partial_y u=0\, ,
$$
hence it parameterizes an $\S^2$ type surface in $\H^3$, having  mean curvature $K$, apart from 
a finite number of branch points.
\end{Lemma}

\proof
Put $\alpha = \frac{1}{2}u_3^{-2}(|{\partial_{x\!}u}|^2 - |{\partial_{y\!}u}|^2)$, $\beta = - u_3^{-2} {\partial_{x\!}u} \cdot {\partial_{y\!}u}$, $\f = \alpha + i \beta$ and
notice that $|\f|\le c_u|\nabla u|^2\in L^\infty(\R^2)$. By direct computation we find
\begin{equation}\label{tech00}
\begin{aligned}
& (\partial_x\alpha- \partial_y\beta) u_3^3 =u_3 {\partial_{x\!}u} \cdot \Delta u - (|{\partial_{x\!}u}|^2- |{\partial_{y\!}u}|^2)\partial_{x\!}u_3 - 2 ({\partial_{x\!}u} \cdot {\partial_{y\!}u})\partial_{y\!}u_3\, , \\
&(\partial_y\alpha + \partial_x\beta) u_3^3 = - u_3 {\partial_{y\!}u} \cdot \Delta u - (|{\partial_{x\!}u}|^2- |{\partial_{y\!}u}|^2)\partial_{y\!}u_3 + 2 ({\partial_{x\!}u} \cdot {\partial_{y\!}u})\partial_{x\!}u_3\, . 
\end{aligned}
\end{equation}
Since $u$ solves \eqref{eq:bollehyp},  it holds that
\begin{equation}\label{tech001}
\begin{aligned}
&u_3 {\partial_{x\!}u} \cdot \Delta u = 2 {G}(\nabla u) \cdot {\partial_{x\!}u} = 2({\partial_{x\!}u} \cdot {\partial_{y\!}u}) \partial_{y\!}u_3 + (|{\partial_{x\!}u}|^2- |{\partial_{y\!}u}|^2)\partial_{x\!}u_3\, ,\\
&u_3 {\partial_{y\!}u} \cdot \Delta u = 2 {G}(\nabla u) \cdot {\partial_{y\!}u} = 2({\partial_{x\!}u} \cdot {\partial_{y\!}u})\partial_{x\!}u_3 - (|{\partial_{x\!}u}|^2- |{\partial_{y\!}u}|^2)\partial_{y\!}u_3\, .
\end{aligned}
\end{equation}
Putting together \eqref{tech00} and \eqref{tech001} we obtain $\partial_x\alpha- \partial_y\beta = \partial_y\alpha + \partial_x\beta = 0$, namely, $\f$ is an holomorphic function. Since $\f$ is bounded and vanishes at infinity
then $\f \equiv 0$ on $\R^2$, hence $u$ is conformal.

The last conclusion follows from Proposition 2.4 and Example 2.5(4) in \cite{GHR}.\QED

\begin{Remark}
\label{R:unbounded}
Here we take $K\equiv k$ constant and point out two simple facts about the energy functional $E_{\text{hyp}}$ in (\ref{eq:ene}). 

By (\ref{eq:Etrinv}), the Nehari manifold contains any nonconstant function. Secondly, $E_{\text{hyp}}$ is unbounded from below. In fact, for $t>1$ we have
$$
E_{\text{hyp}}(\omega+te_3)=\frac12\irn(\omega_3+t)^{-2}~\mu^2\!dz+k\irn(\omega_3+t)^{-2}\omega_3~\mu^2\!dz
=4\pi\big(-\frac{kt-1}{t^2-1}+\frac{k}{2}\ln\frac{t+1}{t-1}\big). 
$$
Notice that $\omega+te_3$ approaches  a horosphere as $t\to 1$, and that $\lim\limits_{t\to 1}E_{\text{hyp}}(\omega+te_3)=-\infty$.
\end{Remark}

\begin{Remark}
Differently from the Euclidean case, see for instance \cite{BrCo2}, the geometric and compactness properties of the energy functional
$E$ are far from being understood (also in the case of a constant curvature), and would deserve a careful analysis.
\end{Remark}

We conclude the paper by pointing out a necessary condition for the existence of embedded $K$ bubbles.

Let $K\in{\cal C}^1(\H^3)$ be given, and let $u\in {\cal C}^2(\ovR^2,\H^3)$ be an embedded solution to (\ref{eq:bollehyp}). By Lemma
\ref{L:conformal}, $u$ is a conformal parametrization of the open set $\Omega\subset\R^3_+$, which is the bounded
connected component of $\R^3_+\setminus u(\S^2)$. We can assume that the nowhere vanishing normal vector $\partial_x u\wedge\partial_y u$ is inward
pointing. Since $u$ is a critical point of the energy functional in (\ref{eq:energy_gen}), then for $j=1,2$ we
have that 
$$
0=E'(u)e_j=V'_K(u)e_j=\irn u_3^{-3}K(u)e_j\cdot \partial_x u\wedge\partial_y u~\!dz=
-\int\limits_{\Omega}\div(p_3^{-3}K(p)e_j)~\!dp
$$
by the divergence theorem. Thus
$$
\int\limits_{\Omega}p_3^{-3}\partial_{p_j} K(p)~\!dp= 0. 
$$
In a similar way, from $E'(u)u=0$ and since $\div(p_3^{-3}K(p)p)= p_3^{-3}\nabla K(p)\cdot p$, one gets
$$
\int\limits_{\Omega}p_3^{-3}\nabla K(p)\cdot p~\!dp= 0.
$$
In particular, $\partial_{p_1} K, \partial_{p_2} K$ and the radial derivative of $K$ can not have constant sign in $\Omega$.
We infer the next nonexistence result (see \cite[Proposition 4.1]{CM} for the Euclidean case).

\begin{Theorem}
\label{T:nonexistence1}
Assume that $K\in {\cal C}^{1}(\H^{3})$ satisfies one of the following conditions,
\begin{itemize}
\item[$i)$] $K(p)=f(\nu\cdot p)$ for some direction $\nu$ orthogonal $e_3$, where $f$ is strictly monotone;
\item[$ii)$] $K(p)=f(|p|)$, where $f$ is strictly monotone.
\end{itemize}
Then (\ref{eq:bollehyp}) has no embedded solution $u\in {\cal C}^2(\ovR^2,\H^3)$.
\end{Theorem}

\bigskip
\noindent
{\bf Acknowledgements}. 
This work is partially supported by PRID-DMIF Projects PRIDEN and VAPROGE, 
 Universit\`a di Udine.

\bigskip

{\footnotesize 
\noindent
{\sc Gabriele Cora},  
Dipartimento di Scienze Matematiche, Informatiche e Fisiche, Universit\`a di Udine
\\ Email: {gabriele.cora@uniud.it}.

\medskip
\noindent
{\sc Roberta Musina}, 
Dipartimento di Scienze Matematiche, Informatiche e Fisiche, Universit\`a di Udine
\\ Email: {roberta.musina@uniud.it}.}
\end{document}